\documentclass[12pt]{amsart}
\usepackage[utf8]{inputenc}
\usepackage{amsmath,amsfonts,amsthm, amssymb, xcolor}
\usepackage{mathtools}
\usepackage{graphicx}
\graphicspath{ {./images/} }
\usepackage{tikz}
\usetikzlibrary{arrows}
\usepackage{float}
\usepackage{caption}
\usepackage{subcaption}
\usepackage{paralist}
\usepackage[pagebackref]{hyperref}
\usepackage[nameinlink]{cleveref}

\usepackage{algorithm}
\usepackage[noend]{algpseudocode}

\usepackage{relsize}

 \newcommand{\freeprod}{\mathop{\Huge \mathlarger{\mathlarger{\ast}}}}

\newtheorem{lemma}{Lemma}[section]
\newtheorem{theorem}{Theorem}[section]
\newtheorem{proposition}{Proposition}[section]
\theoremstyle{definition}

\newtheorem{definition}{Definition}[section]
\newtheorem{corollary}{Corollary}[section]
\newtheorem{remark}{Remark}[section]

\newtheorem{example}{Example}[section]
\newtheorem{construction}{Construction}[section]

\hypersetup{
    colorlinks=true,
    linkcolor=blue,
    filecolor=magenta,      
    urlcolor=cyan,
}

\usepackage[verbose=true,letterpaper]{geometry}
\AtBeginDocument{
  \newgeometry{
    textheight=9in,
    textwidth=6.5in,
    top=1in,
    headheight=12pt,
    headsep=25pt,
    footskip=30pt
  }
}

\title{On cycles and merge trees}
\author{Julian Br{\"u}ggemann and Nicholas A. Scoville}

\address[Julian Br{\"u}ggemann]{Max Planck Institute for Mathematics, Bonn, Germany}
\email{brueggemann@mpim-bonn.mpg.de}

\address[Nicholas A. Scoville]{Department of Mathematics and Computer Science, Ursinus College, Collegeville PA 19426}
\email{nscoville@ursinus.edu}

\date{\today}

\begin{document}
\maketitle
\begin{abstract}
In this paper, we extend the notion of a merge tree to that of a generalized merge tree, a merge tree that includes 1-dimensional cycle birth information.  Given a discrete Morse function on a $1$-dimensional CW complex, i.e., a multigraph, we construct the induced generalized merge tree.  We give several notions of equivalence of discrete Morse functions based on the induced generalized merge tree and how these notions relate to one another. As a consequence, we obtain a complete solution to the inverse problem between discrete Morse functions on $1$-dimensional CW complexes and generalized merge trees. After characterizing which generalized merge trees can be induced by a discrete Morse function on a simple graph, we give an algorithm based on the induced generalized merge tree of a discrete Morse function $f\colon X \to \mathbb{R}$ that cancels the critical cells of $f$ and replaces it with an optimal discrete Morse function.

\end{abstract}
\tableofcontents

\section{Introduction}
Let $X$ be a simplicial complex along with a sequence of subcomplexes $\emptyset=X_0\subseteq X_1 \subseteq \cdots \subseteq X_n=X$  known as a filtration. In the burgeoning field of topological data analysis, a filtration of the standard $n$ simplex is often induced by a point cloud of $n+1$ points based on some increasing parameter related to disctance between points.  Geometrical and topological features of $X$ are then estimated by studying the persistence\footnote{Persistence refers to the range of parameters for which a non-trivial topological feature exists.} of certain topological features \cite{TopPersistence2020}. When the topological feature in question is the number of connected components, the persistence over the lifetime of the filtration is given by birth and death information and is summarized  in a degree 0 barcode or persistence diagram \cite{Oudot2015, Gunnar2022}. If one wishes to not only determine birth and death information from the filtration but also how the components are evolving, that is, which components are merging with which, one associates a merge tree tree to the filtration.  Because the merge tree carries with it this extra information, merge trees are a rich topic of study in both the theoretical and computational settings \cite{2Curry2022, Curry2018, morozov2013interleaving, gasparovic38intrinsic, Curry2022}. Merge trees were introduced as an approximation to contour trees, a special case of Reeb graphs, in the context of visualization \cite{Contour}.

One way to induce a filtration on $X$ is with a discrete Morse function (dMf) \cite{Forman98, Forman2002}. Such a function $f$ induces a filtration by considering subcomplexes associated to each critical value of $f$. The induced merge tree of a dMf on a tree, or 1-dimensional acyclic complex, was introduced in \cite{JohnsonScoville2022}.  There the authors showed that a certain class of merge trees could be realized as the induced merge tree of a star graph.  The authors went on to conjecture that any merge tree could be the induced merge tree of a certain dMf on a path.  This conjecture was recently proved in \cite{brueggemann2021merge}.

The goal of this paper is to extend the theory of merge trees and discrete Morse theory to include cycles. More specifically, given any 1-dimensional CW complex (i.e.~a graph with or without multiedges) equipped with a dMf, we define a generalized induced Morse labeled merge tree (\Cref{iMLT}) associated to this dMf.  The generalized induced Morse labeled merge tree keeps track of not only component birth, death, and merge information but also cycle birth information via a node with a single child. After defining some basic properties, we introduce an equivalence relation on connected graphs called component-merge equivalence (CM equivalence,  \Cref{cmequiv}) and show that there is a one-to-one correspondence between the set of CM equivalence classes of dMfs with only critical cells and the set of isomorphism
classes of generalized Morse labeled merge trees (\Cref{DmfongraphvsMltree}).   In addition, we determine when a given generalized merge tree can be realized by an induced Morse function on a graph without multiedges. Unlike the case of merge trees, not all generalized merge trees can be realized.  \Cref{thm: gmt chracterization} gives a simple counting condition for when a generalized merge tree can be realized by a dMf on a simple graph. The proof is constructive and builds off of the merge tree construction in \cite[Theorem 5.9]{brueggemann2021merge}. Finally in \Cref{section cancel critical}, we give an algorithm on merge trees induced by a dMf in order to cancel critical cells of the dMf. The algorithm allows for some options depending on whether one wishes to preserve homeomorphism type of the graph or find an optimal matching. We briefly compare the algorithm to some similar algorithms from the literature \cite{LEWINER2003221,RandScoville2020}. \par
In the last section, we consider possible future directions and applications.

\addtocontents{toc}{\protect\setcounter{tocdepth}{3}}
	\subsection*{Acknowledgements}\mbox{}\par
	The first author would like to thank Max Planck Institute for Mathematics for the great scientific environment he was part of during his stay. The second author was supported by an AMS-Simons Research Enhancement Grant for PUI Faculty. Both authors would like to thank the anonymous reviewers for their valuable comments and suggestions. 
 
\section{Preliminaries on DMfs and Merge Trees}\label{section preliminaries}
\addtocontents{toc}{\protect\setcounter{tocdepth}{1}}
We recall and introduce the necessary notions for this work. In this article, we use the term graph for finite abstract multigraphs, possibly with self-loops. That is, graphs in this work may have multiple edges between two given vertices, and they can have self-loops, i.e., edges of the form $(x,x)$. This notion of graph can be geometrically interpreted as 1-dimensional CW complexes. \par
On the other hand, we will use the term regular\footnote{We use the term regular in this fashion because we consider the graphs to be combinatorial models for topological spaces, i.e., 1-dimensional regular CW complexes. This should not be confused with the use of the term regular in graph theory. } if $X$ does not contain a self-loop, and we call $X$ a simple graph if $X$ is regular and there is at most one edge between two given vertices. Simple graphs correspond to 1-dimensional simplicial complexes. Since we consider graphs as geometric objects, we also use geometric terms like cells, simplices, and faces to describe them. For any graph $X$, we use $v(X), e(X)$, and $b_1(X)$ to denote the number of vertices, edges, and cycles of $X$, respectively. If an edge $e=uv$ for vertices $u$ and $v$, we say that $u$ and $v$ are the endpoints of $e$.\par
One key feature of this work is that, as usual in works related to topological data analysis, the involved filtrations are considered as part of the data of the space under investigation, rather than just a tool to analyze a space. In this work, the filtrations are given by one of the most central notions of the article, namely that of a discrete Morse function.

\begin{definition}[{\cite[Definition 2.2]{JohnsonScoville2022}/\cite[Section 2.1]{Bene}}] \label{Def:DiscreteMorseFunction}
Let $X$ be a finite graph, not necessarily connected. A function $f\colon X \to \mathbb{R}$ is a \emph{discrete Morse function (dMf)} if it fulfills:
\begin{itemize}
	       \item[\textbf{Monotonicity: }] For cells $\sigma \subset \tau$ we have $f(\sigma)\leq f(\tau)$.
           \item[\textbf{Semi-injectivity: }] $\lvert f^{-1}(\{z\})\rvert \leq 2$ for all $z\in \mathbb{R}$.
           \item[\textbf{genericity: }] For cells $\sigma, \tau \in X$, if $f(\sigma)=f(\tau)$, then either $\sigma \subseteq \tau$ or $\tau \subseteq \sigma$ holds.
           \end{itemize}
% called \emph{weakly increasing} if $f(v)\leq f(e)$ whenever any vertex $v$ is a face of an edge $e$.  A  $f\colon X \to \mathbb{R}$ is a weakly increasing function which is at most 2--1 and satisfies the property that if $f(v)=f(e)$, then $v$ is incident with $e$. 
A cell $\sigma$ of $X$ is  \emph{critical} if $\sigma$ is the unique preimage of $f(\sigma)$. Otherwise, $\sigma$ is called \emph{matched}. If $X$ is not regular, we additionally require that all self-loops are critical with respect to $f$.\par
For any $a\in \mathbb{R}$, the \emph{sublevel subcomplex} of $X$ at $a$ is $X_{a}=\{\sigma\in X: f(\sigma)\leq a\}$. The connected component containing $\sigma \in X$ is denoted $X[\sigma]$. We use the notation $X_{a-\varepsilon}$ to denote the sublevel subcomplex of $X$ immediately preceding $a$, i.e., $X_{a-\varepsilon}:=\{\sigma: f(\sigma)<a\}.$
\end{definition}

\begin{example}\label{ex: first dmf}  Define $f\colon X\to \mathbb{R}$ by
$$
\begin{tikzpicture}[scale=2]

\node[inner sep=2pt, circle] (0) at (0,0) [draw] {};
\node[inner sep=2pt, circle] (2) at (1,0) [draw] {};
\node[inner sep=2pt, circle] (5) at (0,1) [draw] {};
\node[inner sep=2pt, circle] (1) at (1,1) [draw] {};
\node[inner sep=2pt, circle] (4) at (2,1) [draw] {};
\node[inner sep=2pt, circle] (9) at (.5,2) [draw] {};

\path[style=semithick] (0) edge node[anchor=north]{{$8$}}(2);
\path[style=semithick] (0) edge node[anchor=east]{{$10$}}(5);
\path[style=semithick] (5) edge node[anchor=north]{{$6$}}(2);
\path[style=semithick] (2) edge node[anchor=west]{{$3$}}(1);
\path[style=semithick] (4) edge node[anchor=south]{{$11$}}(1);
\path[style=semithick] (5) edge node[anchor=south]{{$7$}}(1);
\path[style=semithick] (9) edge node[anchor=west]{{$12$}}(1);
\path[style=semithick] (5) edge node[anchor=east]{{$13$}}(9);

\node[anchor = north]  at (0) {{$0$}};
\node[anchor = north]  at (2) {{$2$}};
\node[anchor = east]  at (5) {{$5$}};
\node[anchor = north west]  at (1) {{$1$}};
\node[anchor = south]  at (4) {{$4$}};
\node[anchor = south]  at (9) {{$9$}};
\end{tikzpicture}
$$

\noindent Then $f$ is a dMf with each value critical. The sublevel complex $X_7$ is given by

$$
\begin{tikzpicture}[scale=2]

\node[inner sep=2pt, circle] (0) at (0,0) [draw] {};
\node[inner sep=2pt, circle] (2) at (1,0) [draw] {};
\node[inner sep=2pt, circle] (5) at (0,1) [draw] {};
\node[inner sep=2pt, circle] (1) at (1,1) [draw] {};
\node[inner sep=2pt, circle] (4) at (2,1) [draw] {};

\path[style=semithick] (5) edge node[anchor=north]{{$6$}}(2);
\path[style=semithick] (2) edge node[anchor=west]{{$3$}}(1);
\path[style=semithick] (5) edge node[anchor=south]{{$7$}}(1);

\node[anchor = north]  at (0) {{$0$}};
\node[anchor = north]  at (2) {{$2$}};
\node[anchor = east]  at (5) {{$5$}};
\node[anchor = north west]  at (1) {{$1$}};
\node[anchor = south]  at (4) {{$4$}};
\end{tikzpicture}
$$
where $X_7[f^{-1}(1)]=X_7[f^{-1}(6)]$ is the 3-cycle while $X_7[f^{-1}(4)]$ is simply the isolated vertex labeled $4.$
\end{example}

\begin{remark}
    This definition of dMfs, due to B. Benedetti, is not equivalent to the more general definition originally given by Forman \cite{Forman98}, but stronger. In particular, every dMf in the sense of \Cref{Def:DiscreteMorseFunction} is also a dMf in Forman's sense. The definition presented above is generic, i.e., any dMf in the sense of Forman can be modified to fulfill the definition above without changing the induced acyclic matching. Nonetheless, adjusting a dMf in Forman's sense to become a dMf in the sense of \Cref{Def:DiscreteMorseFunction} will in general change the associated filtration of the complex at hand. The definition stated above has the advantage that critical cells are distinguished by their critical values and at each level, at most either one critical cell or one pair of matched cells is added to the sublevel complex. We choose to use the term ``dMf'' to follow Benedetti's notation. This should not lead to confusion because all dMfs in this work satisfy \Cref{Def:DiscreteMorseFunction}.\par
    The condition that self-loops must be critical is a standard approach to dMfs on CW complexes that fail to be regular. In detail, one requires non-regular faces to remain unmatched. In the 1-dimensional case, self-loops are the only way to break regularity.  
\end{remark}

While discrete Morse theory provides a well-developed framework for filtered spaces, we apply this framework to one specific kind of topological information, namely, the development of connected components throughout the filtration. The development of connected components is summarized in the merge tree. 

\begin{definition}[{\cite[Definition 3.1]{JohnsonScoville2022}}]\label{Def:MergeTree}
A rooted tree is called \textit{binary} if it is a rooted tree where each vertex has at most two children. A binary tree is called \textit{full} if each node has either 0 or 2 children. A binary tree is called \textit{chiral}\footnote{For full binary trees, the notion of chirality is equivalent to the notion of ordered trees, i.e., trees that for each node have a specified total order on the set children of that node. We define chirality in a more general fashion here in order to use the same term for the definition of general merge trees.} if each vertex is equipped with a label of either L or R with the children of a vertex having one label L and the other R. A \textit{merge tree} is a chiral full binary tree.\par
For a node $p$ of a rooted tree, we denote by $T(p)$ the \textit{rooted subtree with root} $p$, that is, the subtree that consists of $p$ and all of $p$'s descendants. 
\end{definition}

\begin{remark}\label{Rem:GeometricMergeTreeVSOurMergeTree}
    Merge trees were introduced in \cite{Contour}\footnote{In that work, merge trees are called join trees.} to compute contour trees efficiently. Although the original introduction of merge trees was of a combinatorial nature, it is also common in the literature to perceive merge trees of filtration maps $f\colon X \rightarrow \mathbb{R}$ in a more geometric way, that is, as the quotient space $X/\sim$, where $\sim$ is defined by $x\sim y$ if and only if $f(x)=f(y)$ and $x$ and $y$ are in the same connected component of the sublevel set $f^{-1}(-\infty, f(x)]$. If one assumes the filtration function $f$ to be generic\footnote{In this work, this sense of genericity is ensured by the genericity and semi-injectivity properties in \Cref{Def:DiscreteMorseFunction}.} in the sense that at each time at most two connected components merge, the induced merge tree will be binary. From that point of view one obtains a merge tree in the sense of \Cref{Def:MergeTree} in the following way. Consider new connected components in the filtration as leaves and consider points where connected components merge as inner nodes. This way the merge tree in the geometric sense becomes a full binary tree. The chirality in \Cref{Def:MergeTree} is an additional structure that is somewhat motivated by the Elder rule: connected components that were created earlier should persist longer. This is reflected in the way that any (generalized) merge tree constructed according to \Cref{iMLT} has the property that bars in the induced 0-barcode correspond to maximal paths in the induced merge tree that only go through nodes of the same chirality.\par
    In that sense, it would be more adequate to refer to merge trees in the sense of \Cref{Def:MergeTree} as chiral generic merge trees but we decide against that in the interest or brevity because all merge trees in this work are chiral and generic. 
\end{remark}

Our main object of study is given in \Cref{def: generalized merge tree}, that of a generalized merge tree.  

\begin{definition}\label{def: generalized merge tree}
   A \textit{generalized merge tree} $T$ is a chiral binary tree $T$ such that each leaf has a sibling, and inner nodes without a sibling have the same chirality as their parent node. By convention, we say that the root always has chirality L. Furthermore, the root is never regarded as a leaf, even if it only has one child node. A \textit{generalized merge subtree} $T'$ of a generalized merge tree $T$ is a rooted subtree of a node $p\in T$, such that each node of $T'$ has the same chirality as in $T$. 
\end{definition}

A generalized merge tree generalizes the notion of a merge tree in the following sense:  while a merge tree keeps track of component information, the generalized merge tree will also keep track of 1-dimensional cycle information.  A cycle is represented by a child vertex with no sibling.  See \Cref{ex: generalized Morse merge tree}.

For nodes $c$ of generalized merge trees we use the notation $\textcolor{red}{c_{l}}/\textcolor{blue}{c_r}$ for the $\textcolor{red}{left}/\textcolor{blue}{right}$ child node of $c$.

\begin{remark}
Generalized merge trees may have nodes without siblings. Thus they are allowed  to break the condition of being full since some nodes might have only one child node. This breaks the property of being full in a different way than the usual notion of merge tree does: non-generic\footnote{In the sense as in \Cref{Rem:GeometricMergeTreeVSOurMergeTree}.} filtration functions induce merge trees where each node might have two or more children but single children will never occur. 
% The notion of chirality of only children does not really deserve the name chirality because there is only one total ordering on a set with one element. 
We impose the condition that an only child has the same chirality as its parent node for technical reasons. We need this convention so \Cref{idmf} and \Cref{iMLT} make part 2 of \Cref{DmfongraphvsMltree} work.  \par
\end{remark}

\begin{example}\label{ex: chrial binary nerge tree}
The tree $T$ below is a generalized merge tree:

$$
\begin{tikzpicture}[scale=1]

\node[inner sep=2pt, circle] (13) at (0,0) [draw] {};
\node[inner sep=2pt, circle] (12) at (-1,1) [draw] {};
\node[inner sep=2pt, circle] (9) at (0,2) [draw] {};
\node[inner sep=2pt, circle] (11) at (-2,2) [draw] {};
\node[inner sep=2pt, circle] (4) at (-1,3) [draw] {};
\node[inner sep=2pt, circle] (10) at (-3,3) [draw] {};
\node[inner sep=2pt, circle] (8) at (-4,4) [draw] {};
\node[inner sep=2pt, circle] (0) at (-5,5) [draw] {};
\node[inner sep=2pt, circle] (7) at (-3,5) [draw] {};
\node[inner sep=2pt, circle] (6) at (-2,6) [draw] {};
\node[inner sep=2pt, circle] (5) at (-3,7) [draw] {};
\node[inner sep=2pt, circle] (3) at (-1,7) [draw] {};
\node[inner sep=2pt, circle] (1) at (0,8) [draw] {};
\node[inner sep=2pt, circle] (2) at (-2,8) [draw] {};

\path[style=semithick] (12) edge node[anchor=north]{{}}(13);
\path[style=semithick] (12) edge node[anchor=west]{{}}(9);
\path[style=semithick] (12) edge node[anchor=south]{{}}(11);
\path[style=semithick] (11) edge node[anchor=south]{{}}(4);
\path[style=semithick] (11) edge node[anchor=south]{{}}(10);
\path[style=semithick] (8) edge node[anchor=south]{{}}(10);
\path[style=semithick] (8) edge node[anchor=south]{{}}(0);
\path[style=semithick] (8) edge node[anchor=south]{{}}(7);
\path[style=semithick] (7) edge node[anchor=south]{{}}(6);
\path[style=semithick] (6) edge node[anchor=south]{{}}(5);
\path[style=semithick] (6) edge node[anchor=south]{{}}(3);
\path[style=semithick] (3) edge node[anchor=south]{{}}(1);
\path[style=semithick] (3) edge node[anchor=south]{{}}(2);

\node[anchor = north]  at (13) {L};
\node[anchor = north]  at (12) {L};
\node[anchor = north]  at (11) {L};
\node[anchor = north]  at (10) {L};
\node[anchor = north]  at (8) {L};
\node[anchor = north]  at (0) {L};
\node[anchor = south]  at (4) {R};
\node[anchor = south]  at (9) {R};
\node[anchor = west]  at (7) {R};
\node[anchor = west]  at (6) {R};
\node[anchor = west]  at (3) {R};
\node[anchor = west]  at (1) {R};
\node[anchor = east]  at (2) {L};
\node[anchor = east]  at (5) {L};
\end{tikzpicture}
$$
Note that each vertex has at most one child and that inner nodes without a sibling have the same chirality (label of L or R) as their parent.  This generalized merge tree is said to have 6 leaves, as the root node at the bottom (labeled L) is not considered a leaf by convention. \par
Any generalized merge subtree $T'$ of $T$ is obtained by picking a node $p$ of $T$ and taking all descendants of $p$ together with their chirality. For example, the following is a generalized merge subtree of $T$:
$$
\begin{tikzpicture}
   \node[inner sep=2pt, circle] (7) at (-3,5) [draw] {};
\node[inner sep=2pt, circle] (6) at (-2,6) [draw] {};
\node[inner sep=2pt, circle] (5) at (-3,7) [draw] {};
\node[inner sep=2pt, circle] (3) at (-1,7) [draw] {};
\node[inner sep=2pt, circle] (1) at (0,8) [draw] {};
\node[inner sep=2pt, circle] (2) at (-2,8) [draw] {};

\path[style=semithick] (7) edge node[anchor=south]{{}}(6);
\path[style=semithick] (6) edge node[anchor=south]{{}}(5);
\path[style=semithick] (6) edge node[anchor=south]{{}}(3);
\path[style=semithick] (3) edge node[anchor=south]{{}}(1);
\path[style=semithick] (3) edge node[anchor=south]{{}}(2);

\node[anchor = west]  at (7) {R};
\node[anchor = west]  at (6) {R};
\node[anchor = west]  at (3) {R};
\node[anchor = west]  at (1) {R};
\node[anchor = east]  at (2) {L};
\node[anchor = east]  at (5) {L};
\end{tikzpicture}
$$
One needs to be careful that generalized merge subtrees are not necessarily generalized merge trees because their roots might have chirality $R$.
    
\end{example}
% \begin{definition}[{\cite[Definition 2.17]{brueggemann2021merge}}]\label{morseorder}
%      Let $T$ be a generalized merge tree. We call a total order $\leq$ on the nodes of $T$ a \textit{Morse order} if it fulfills the following two properties for all generalized merge subtrees $T'$ of $T$:
%      \begin{enumerate}
%         \item The restriction $\leq_{\lvert T'}$ attains its maximum on the root $p$ of $T'$.
%         \item The minimum of $\leq_{\lvert T'}$ has the same chirality as $p$.
%     \end{enumerate}
%     We call any pair $(T,\leq)$ of a generalized merge tree $T$ together with a Morse order $\leq$ a \emph{generalized Morse ordered merge tree (gMO tree)}.
%     \end{definition}
    
\begin{definition}[{\cite[Definition 2.19]{brueggemann2021merge}}]\label{morselabel}
 We call a generalized merge tree $(T,\lambda)$ with an injective map $\lambda \colon T \rightarrow \mathbb{R}$ such that
 \begin{enumerate}
        \item for all generalized merge subtrees $T'$ of $T$ the restriction $\lambda_{\lvert T'}$ attains its maximum on the root $p$ of $T'$.
        \item for all generalized merge subtrees $T'$ of $T$ the node on that $\lambda_{\lvert T'}$ attains its minimum has the same chirality as $p$.
    \end{enumerate}
 % $\lambda$ induces
 % \footnote{By pulling back the total order of the reals.}  a Morse order on $T$ 
 a \textit{generalized Morse labeled merge tree (gML tree)}. Any such map $\lambda$ is called a \textit{Morse labeling} on $T$.\par
\end{definition}
\begin{remark}\label{morseord}
	Property 2 of \Cref{morselabel} is equivalent to either of the following:
	\begin{itemize}
		\item For any generalized merge subtrees $T'$ with root $p$ of $T$, the restriction $\leq_{\lvert T'}$ attains its minimum on the subtree with root $\textcolor{red}{p_{l}}/\textcolor{blue}{p_r}$ if \textcolor{red}{L}/\textcolor{blue}{R} is the chirality of the root $p$ of $T'$.
		\item For any generalized merge subtrees $T'$ with root $p$ of $T$, all nodes on the shortest path between $p$ and the minimum of $\leq_{\lvert T'}$ have the same chirality as $p$. 
	\end{itemize}
	The equivalence can be proved by an inductive argument over all nodes of the shortest path between $p$ and the minimum.\par
    We refer to the final step of \Cref{ex: generalized Morse merge tree} for an example of a 
    % gMO tree as well as a 
    gML tree. For any generalized merge tree, there is always at least one Morse labeling, usually several, see e.g.~\Cref{scmo} or with minor modifications \cite[Definition 3.3]{brueggemann2021merge}.
\end{remark}
\begin{definition}
 Let $(T,\lambda)$ and $(T',\lambda')$ be gML trees. An \textit{order equivalence} $(\varphi,\psi)\colon(T,\lambda)\rightarrow (T',\lambda')$ of gML trees is a pair of maps consisting of an isomorphism of the underlying generalized merge trees $\varphi\colon T\rightarrow T'$ and a bijection $\psi\colon \mathbb{R}\rightarrow\mathbb{R}$ such that the restriction $\psi_{\lvert \operatorname{im}(\lambda)}\colon \operatorname{im}(\lambda)\rightarrow \operatorname{im}(\lambda')$ is order preserving.
\end{definition}

% \begin{proposition}\label{MorselabelvsMorseorder}
% Let $gMOT$ be the set of generalized Morse ordered merge trees and let $gMLT$ be the set of generalized Morse labeled merge trees. Then (i) taking the Morse order induced by a Morse labeling and (ii) using a Morse order with labels $\{0,\cdots, |V(T)|-1\}$ to induce a Morse labeling, define inverse bijections
% $$
% \begin{tikzpicture}
% \node (0) at (0,0) {$iML \colon gMOT/_{\cong}$}; 
% \node (1) at (4,0) {$gMLT/_{\sim}\colon iMO$};
% \draw[<->] (0) -- (1);
% \end{tikzpicture}
% $$
% where $\sim$ denotes order equivalence.
% \end{proposition}
% \begin{proof}
% The proof is analogous to \cite[Proposition 3.26]{brueggemann2021merge}.
% \end{proof}
\begin{definition}
    Let $(X,f)$ be a dMf on a graph. We call a critical edge $\sigma \in X$ a \textit{closing edge} if there is a subdivision $\Delta$ of $S^1$ contained as a subcomplex $\Delta \subseteq X$ which contains $\sigma$ such that $f(\sigma)$ is the maximum of $f$ on $\Delta$.\par
    We define $C(X,f)\coloneqq \{ c \in X \vert c \text{ is closing}\}$ to be the \textit{set of closing edges} of $(X,f)$ and $(\Bar{X},\Bar{f})\coloneqq (X \setminus C(X,f),f_{\lvert X \setminus C(X,f)})$ to be the \textit{spanning tree (or spanning forest, if $X$ is not connected) induced by $f$} of $X$. 
\end{definition}
\begin{remark}
    In the previous definition, the subdivison of $S^1$ that any closing edge $\sigma$ must be part of does not need to be unique. Nonetheless, the removal of $\sigma$ would lead to the reduction of the first Betti number by one. Moreover, the notion of closing edges is well-defined because the edge $\sigma$ being closing implies that it is the unique maximal edge of all subdivisions of $S^1$ in $X_{f(\sigma)}[\sigma]$ that contain $\sigma$.\par
    Furthermore, it is immediate that $(\Bar{X},\Bar{f})\coloneqq (X \setminus C(X,f),f_{\lvert X \setminus C(X,f)})$ is a dMf on a tree. It is also immediate that self-loops are always closing edges. 
\end{remark}

\begin{construction}[Induced gML Merge Tree\footnote{This construction generalizes \cite[Theorem 3.5]{JohnsonScoville2022}.}]\label{iMLT}
Let $f\colon X \rightarrow \mathbb{R}$ be a dMf on a connected graph $X$. If $f$ has a single critical vertex $v$ with $f(v)=r$ and no critical edges, then the induced gML merge tree is a single node $c_{\lambda_v}$  with label $\lambda_f(c_{\lambda_v})=r$ and chirality L.  Otherwise, let $\sigma_n>\sigma_{n-1}> \dots > \sigma_1>\sigma_0$ be the critical edges of $(X,f)$ ordered by their values under $f$. The gML merge tree induced by $(X,f)$, denoted $M(X,f)$ with labeling $\lambda_f$, is constructed inductively, inducing over the decreasing order of the critical edges.\\

Start by constructing a root node called $M(\sigma_n)$, labeled $\lambda_f(M(\sigma_n)):=f(\sigma_n)$, and left chirality. 

Now begin the induction over the decreasing order of the critical edges starting from $n,\ldots, 0$. For  $\sigma_i$ a critical edge with endpoints $u$ and $v$, the addition of $\sigma_i$ at level subcomplex $X_{f(\sigma_i)}$ either creates a cycle or connects two components.  Formally,
\begin{enumerate}
    \item The critical edge $\sigma_i$ is closing.\footnote{This condition is equivalent to $b_1(X_{f(\sigma_i)} \setminus \sigma_i)=b_1(X_{f(\sigma_i)})-1$, which in turn is equivalent to $b_0(X_{f(\sigma_i)}\setminus \sigma_i)=b_0(X_{f(\sigma_i)})$.}
    % $\Leftrightarrow$ $b_1(X_{f(\sigma_i)} \setminus \sigma_i)=b_1(X_{f(\sigma_i)})-1$ $\Leftrightarrow$ $b_0(X_{f(\sigma_i)}\setminus \sigma_i)=b_0(X_{f(\sigma_i)})$,  
    \item The critical edge $\sigma_i$ is not closing.\footnote{This condition is equivalent to $b_1(X_{f(\sigma_i)}\setminus \sigma_i)=b_1(X_{f(\sigma_i)})$, which in turn is equivalent to $b_0(X_{f(\sigma_i)}\setminus \sigma_i)=b_0(X_{f(\sigma_i)}) +1$.}
    % $\Leftrightarrow$ $b_1(X_{f(\sigma_i)}\setminus \sigma_i)=b_1(X_{f(\sigma_i)})$ $\Leftrightarrow$ $b_0(X_{f(\sigma_i)}\setminus \sigma_i)=b_0(X_{f(\sigma_i)}) +1$,
\end{enumerate}
If $\sigma_i$ satisfies (1), we construct a child node $c$ of $M(\sigma_i)$ with label $\lambda_f \coloneqq \max \{f(\sigma)\vert \sigma \in X_{f(\sigma_i)-\varepsilon}, \ \sigma \text{ is critical} \}$ and the same chirality as $M(\sigma_i)$. The node $c$ then corresponds to the edge of $X$ labeled $\lambda$.

\noindent If $\sigma_i$ satisfies (2), we construct two child nodes $c_{\lambda_v}$ and $c_{\lambda_w}$ of $M(\sigma_i)$. Define $\lambda_v\coloneqq$ $\max\{f(\sigma) \rvert \sigma \in X_{f(\sigma_i)-\varepsilon}[v], \sigma \text{ critical}\}$ and $\lambda_w\coloneqq\max\{f(\sigma) \vert \sigma \in X_{f(\sigma_i)-\varepsilon}[w], \sigma \text{ is critical.}\}$.
    Then label the new nodes $\lambda_f(c_{\lambda_v})\coloneqq\lambda_v$ and $\lambda_f(c_{\lambda_w})\coloneqq\lambda_w$. If $\min\{f(\sigma) \vert \sigma \in X_{f(\sigma_i)-\varepsilon}[v]\}<\min\{f(\sigma) \vert \sigma \in X_{f(\sigma_i)-\varepsilon}[w]\}$, we assign $c_{\lambda_v}$ the same chirality (L or R) as $c_{\sigma_i}$ and give $c_{\lambda_w}$ the opposite chirality.\par
    
    Continue the induction over all of the critical edges of $X$ to obtain the Morse labeled merge tree $M(X,f)$ induced by $f$ along with labeling $\lambda_f$.

\end{construction}

\begin{remark}\label{Rem:CorrespondenceBetweenNodesAndVerticesAndEdges}
It is important to note that in the inductive step after creating the child(ren) of $M(\sigma_i)$, if critical edge $\sigma_{i-1}$ exists, then $M(\sigma_{i-1})$ is a node on the current constructed generalized merge tree (it may be a child of a node other than $M(\sigma_i)$). Leaf nodes, for which no such critical edge exists, correspond to critical vertices of $f$. That is, after the final inductive step all the leaf nodes correspond to the critical vertices of $f$ and all the inner nodes correspond to the critical edges of $f$. 
\end{remark}

\begin{example}\label{ex: generalized Morse merge tree}
We will construct the induced Morse labeled merge tree of \Cref{ex: first dmf}.  Since the dMf in this example is injective, we will name each vertex or edge by $\sigma_i$ where $f(\sigma)=i$.  The first step of \Cref{iMLT} is to list the critical edges in increasing order:
$$
\sigma_{13}>\sigma_{12}>\sigma_{11}>\sigma_{10}>\sigma_{8}>\sigma_{7}>\sigma_{6}.
$$
Here we index by the value under the dMf as opposed to the integers $6, \ldots, 0$ but it does not matter.  For the base case, we create a node called $M(\sigma_{13})$ with value $\lambda_f(M(\sigma_{13}))=f(\sigma_{13})=13$ with chirality L (by definition); that is, we begin with

$$
\begin{tikzpicture}[scale=1]

\node[inner sep=2pt, circle] (13) at (0,0) [draw] {};

\node[anchor = north]  at (13) {13L};

\end{tikzpicture}
$$

Moving on from the base case, $\sigma_{13}$ creates a cycle, i.e., it is a closing edge so that it has a single child $M(\sigma_{12})$ with label $\lambda_f=\max\{f(\sigma): \sigma \in X_{13-\epsilon}, \sigma \text{ critical }\}=12$ and chirality that of $M(\sigma_{13})$ which is L.  This yields

$$
\begin{tikzpicture}[scale=1]

\node[inner sep=2pt, circle] (13) at (0,0) [draw] {};
\node[inner sep=2pt, circle] (12) at (-1,1) [draw] {};

\path[style=semithick] (12) edge node[anchor=north]{{}}(13);

\node[anchor = north]  at (13) {13L};
\node[anchor = east]  at (12) {12L};

\end{tikzpicture}
$$
Now since $M(\sigma_{12})$ is not a closing edge, we construct two child nodes $c_{\lambda_{\sigma_9}}$ and $c_{\lambda_{\sigma_1}}$ of $M(\sigma_{12})$.  There values are then computed as
$$
\lambda_{\sigma_9}=\max\{f(\sigma): \sigma\in X_{12-\epsilon}[\sigma_9], \sigma \text{ critical } \}=9
$$
and
$$
\lambda_{\sigma_1}=\max\{f(\sigma): \sigma\in X_{12-\epsilon}[\sigma_1], \sigma \text{ critical } \}=11.
$$
This amounts to determining the largest critical value in the connected component of the vertex in question. Hence the two children of $M(\sigma_{12})$ are labeled 9 and 11.  Finally, $0=\min\{f(\sigma): \sigma\in X_{12-\epsilon}[\sigma_1], \sigma \text{ critical } \}<\min\{f(\sigma): \sigma\in X_{12-\epsilon}[\sigma_9], \sigma \text{ critical } \}=9$ so that $c_{\lambda_{\sigma_1}}$ shares the same chirality as its parent while $c_{\lambda_{\sigma_9}}$ has the opposite chirality.  In sum, we have so far

$$
\begin{tikzpicture}[scale=1]

\node[inner sep=2pt, circle] (13) at (0,0) [draw] {};
\node[inner sep=2pt, circle] (12) at (-1,1) [draw] {};
\node[inner sep=2pt, circle] (9) at (0,2) [draw] {};
\node[inner sep=2pt, circle] (11) at (-2,2) [draw] {};

\path[style=semithick] (12) edge node[anchor=north]{{}}(13);
\path[style=semithick] (12) edge node[anchor=west]{{}}(9);
\path[style=semithick] (12) edge node[anchor=south]{{}}(11);

\node[anchor = north]  at (13) {13L};
\node[anchor = east]  at (12) {12L};
\node[anchor = east]  at (11) {11L};
\node[anchor = north]  at (9) {9R};

\end{tikzpicture}
$$
The induction again continues at $\sigma_{11}$ which is not a closing edge.  The two child nodes of $M(\sigma_{11})$ have values $4$ and $10$ with the nodes given value 4 sharing the same chirality as $M(\sigma_{11})$ so that we have

$$
\begin{tikzpicture}[scale=1]

\node[inner sep=2pt, circle] (13) at (0,0) [draw] {};
\node[inner sep=2pt, circle] (12) at (-1,1) [draw] {};
\node[inner sep=2pt, circle] (9) at (0,2) [draw] {};
\node[inner sep=2pt, circle] (11) at (-2,2) [draw] {};
\node[inner sep=2pt, circle] (4) at (-1,3) [draw] {};
\node[inner sep=2pt, circle] (10) at (-3,3) [draw] {};

\path[style=semithick] (12) edge node[anchor=north]{{}}(13);
\path[style=semithick] (12) edge node[anchor=west]{{}}(9);
\path[style=semithick] (12) edge node[anchor=south]{{}}(11);
\path[style=semithick] (11) edge node[anchor=south]{{}}(4);
\path[style=semithick] (11) edge node[anchor=south]{{}}(10);

\node[anchor = north]  at (13) {13L};
\node[anchor = east]  at (12) {12L};
\node[anchor = east]  at (11) {11L};
\node[anchor = east]  at (10) {10L};
\node[anchor = north]  at (4) {4R};
\node[anchor = north]  at (9) {9R};

\end{tikzpicture}
$$
Continuing in this manner we arrive at the induced Morse labeled merge tree given by

$$
\begin{tikzpicture}[scale=1]

\node[inner sep=2pt, circle] (13) at (0,0) [draw] {};
\node[inner sep=2pt, circle] (12) at (-1,1) [draw] {};
\node[inner sep=2pt, circle] (9) at (0,2) [draw] {};
\node[inner sep=2pt, circle] (11) at (-2,2) [draw] {};
\node[inner sep=2pt, circle] (4) at (-1,3) [draw] {};
\node[inner sep=2pt, circle] (10) at (-3,3) [draw] {};
\node[inner sep=2pt, circle] (8) at (-4,4) [draw] {};
\node[inner sep=2pt, circle] (0) at (-5,5) [draw] {};
\node[inner sep=2pt, circle] (7) at (-3,5) [draw] {};
\node[inner sep=2pt, circle] (6) at (-2,6) [draw] {};
\node[inner sep=2pt, circle] (5) at (-3,7) [draw] {};
\node[inner sep=2pt, circle] (3) at (-1,7) [draw] {};
\node[inner sep=2pt, circle] (1) at (0,8) [draw] {};
\node[inner sep=2pt, circle] (2) at (-2,8) [draw] {};

\path[style=semithick] (12) edge node[anchor=north]{{}}(13);
\path[style=semithick] (12) edge node[anchor=west]{{}}(9);
\path[style=semithick] (12) edge node[anchor=south]{{}}(11);
\path[style=semithick] (11) edge node[anchor=south]{{}}(4);
\path[style=semithick] (11) edge node[anchor=south]{{}}(10);
\path[style=semithick] (8) edge node[anchor=south]{{}}(10);
\path[style=semithick] (8) edge node[anchor=south]{{}}(0);
\path[style=semithick] (8) edge node[anchor=south]{{}}(7);
\path[style=semithick] (7) edge node[anchor=south]{{}}(6);
\path[style=semithick] (6) edge node[anchor=south]{{}}(5);
\path[style=semithick] (6) edge node[anchor=south]{{}}(3);
\path[style=semithick] (3) edge node[anchor=south]{{}}(1);
\path[style=semithick] (3) edge node[anchor=south]{{}}(2);

\node[anchor = north]  at (13) {13L};
\node[anchor = east]  at (12) {12L};
\node[anchor = east]  at (11) {11L};
\node[anchor = east]  at (10) {10L};
\node[anchor = east]  at (8) {8L};
\node[anchor = east]  at (0) {0L};
\node[anchor = west]  at (4) {4R};
\node[anchor = west]  at (9) {9R};
\node[anchor = west]  at (7) {7R};
\node[anchor = west]  at (6) {6R};
\node[anchor = west]  at (3) {3R};
\node[anchor = west]  at (1) {1R};
\node[anchor = east]  at (2) {2L};
\node[anchor = east]  at (5) {5L};
\end{tikzpicture}
$$
which is the same merge tree as in \Cref{ex: chrial binary nerge tree}.

\end{example}

\begin{remark}\label{rem graph and tree correspondence}
The construction of the induced gML tree comes with a bijection $M\colon X \rightarrow V(M(X,f))$ that restricts to bijections between the critical vertices of $X$ and leaves of $M(X,f)$, between the non-closing critical edges of $X$ and parents with two children of $M(X,f)$, and between the closing edges (cycles) of $X$ and parents with one child in $M(X,f)$.\par
Furthermore, the proof that the construction indeed produces a gML tree is completely analogous to \cite[Proposition 2.20]{brueggemann2021merge}, respectively \cite[Theorem 9]{JohnsonScoville2022}.\par
It is also possible to apply the construction to dMfs on non-connected graphs. In that case the algorithm produces a merge forest and one can deal with each connected component separately.
\end{remark}

\begin{lemma}\label{lem: subtree component}
Let $(X,f)$ be a dMf on a graph and let $M(X,f)$ be the induced gML tree. For any critical cell $s\in X$, the generalized merge subtree\footnote{Recall \Cref{iMLT} for the definition of $M(s)$ and \Cref{def: generalized merge tree} for the definition of generalized merge subtrees $T(\_)$.} $T(M(s))$ of $M(X,f)$ is induced by the connected component $X_{f(s)}[s]$ of $s$ in the sublevel complex of level $f(s)$. Moreover, the rooted subtree $T(M(s))$ is isomorphic to $M(X_{f(s)}[s],f_{\lvert X_{f(s)}[s]})$ as merge trees if and only if $M(s)$ has chirality L. If $M(s)$ has chirality R, then $T(M(s))$ is isomorphic to $M(X_{f(s)}[s],f_{\lvert X_{f(s)}[s]})$ as rooted binary trees but the chiralities of all nodes are opposite to the ones of their respective nodes in the other tree. 
\end{lemma}

\begin{proof}
We observe that by \Cref{iMLT} the label of $M(s)$ is $f(s)$ and the chirality of $M(s)$ is decided by the minimum of $f_{\lvert X_{f(s)}[s]}$ in comparison to the minimum of the connected component that $X_{f(s)}[s]$ was divided from at level $f(s)$. It follows inductively by construction that all nodes of the subtree $T(M(s))$ are induced by critical cells of $X_{f(s)}[s]$ because they are constructed by removing critical edges of $X_{f(s)}[s]$.\\
The isomorphism as rooted binary trees is constructed by the same inductive argument. Since the chirality depends on the chirality of the respective parent node, said isomorphism is compatible with the chirality if and only if the root of the rooted subtree $T(M(s))$, namely $M(s)$, has chirality L. This is true because the root of $M(X_{f(s)}[s],f_{\lvert X_{f(s)}[s]})$ by convention always has chirality L.
\end{proof}
\begin{definition}\label{Def:CycleNodesandUnderlyingMLTree}
Let $(T,\lambda)$ be a gML tree. Let $C(T)\subset V(T)$ be the set of nodes that have exactly one child node. We refer to the elements of $C(T)$ as \emph{cycle nodes}. \par
We denote by $(\Bar{T},\lambda)$ the Morse labeled merge tree that is obtained from $(T,\lambda)$ by removing the cycle nodes by connecting their parent nodes directly to their child nodes. We call $(\Bar{T},\lambda)$ the \emph{underlying Morse labeled merge tree} of $(T,\lambda)$.
\end{definition}
% \begin{definition}[{\cite[Definition 3.1]{brueggemann2021merge}}]\label{pw}
% 	Let $T$ be a generalized merge tree of depth $n$ and let $a$ be a node of $T$. The \textit{path word} corresponding to $a$ is a word $a_0a_1\dots a_n \in \{L,R,\_\}^{n+1}$ where $\_$ denotes the empty letter. If $a$ is of depth $k$, the letters $a_0\dots a_k$ are given by the chirality of the nodes belonging to the shortest directed path from the root to $a$. The letters $a_{k+1}\dots a_n$ are then empty.
% \end{definition}
% \begin{remark}\label{maxk}
%     Let $a,b$ be nodes of a generalized merge tree $T$ and let $a_0a_1\dots a_n$ be the path word corresponding to $a$ and $b_0b_1\dots b_n$ be the path word corresponding to $b$. Then the equation $a_0=b_0=L$ always holds because we consider paths that begin at the root. Because $a_0=b_0=L$ and because we consider finite trees, there is always a maximal $k \in \mathbb{N}$ such that $a_i=b_i$ holds for all $i\leq k$. Furthermore, the last non-empty letter of a path word is always the chirality of the considered node. 
% \end{remark}
\begin{definition}[{\cite[Definition 3.20]{brueggemann2021merge}}]\label{idMfonPath}
Let $(T,\lambda)$ be a Morse labeled merge tree. Let $P$ be an oriented path with $i(T)$ 1-simplices, where $i(T)$ is the number of inner nodes of $T$. We consider $P$ as a totally ordered set $(P,\leq)$ by saying that a simplex $\sigma$ is less than a simplex $\tau$ if $\sigma$ is further left than $\tau$, see \cite[Definition 3.15]{brueggemann2021merge} for details. Moreover, we consider $T$ as a totally ordered set\footnote{With a different order than induced by the Morse labeling!} by saying that a node $c$ is less than a node $c'$ if $c$ is left of $c'$ in $T$, see \cite[Definition 3.11]{brueggemann2021merge} for details. We denote the unique order preserving bijection, see \cite[Remark 3.18]{brueggemann2021merge}, between $P$ and $T$ by $\phi\colon P \rightarrow T$. Then the discrete Morse function on $P$ induced by $(T,\lambda)$ is $f_\lambda\coloneqq \lambda \circ \phi$.
\end{definition}
\begin{construction}[Induced dMf]\label{idmf}
Let $(T,\lambda)$ be a gML tree. We obtain a dMf on a graph $f_\lambda\colon X \rightarrow \mathbb{R}$ from $(T,\lambda)$ in two steps as follows:
In a first step, we consider the underlying Morse labeled merge tree $(\Bar{T},\lambda)$ and construct its induced dMf on a path $(P,f_\lambda)$ as in \Cref{idMfonPath}. For the second step, for each node $c$ of $C(T)$ we add an edge parallel to the edge corresponding to $c$'s oldest descendant which has two children to $P$. We denote the graph obtained this way by $X$ and extend the function $f_\lambda\colon P\rightarrow \mathbb{R}$ to $X$ using the values of $\lambda$ on the corresponding nodes. We denote the pair $(X,f_\lambda)$ by $\Phi(T,\lambda)$ and consider that we also obtained a bijection $\phi \colon V(T) \rightarrow \Phi(T,\lambda)$.
\end{construction}
\begin{example}\label{ex:induced dMf}
    We construct the induced discrete Morse function of the gML tree from \Cref{ex: generalized Morse merge tree}. The underlying merge tree is:
    $$
\begin{tikzpicture}[scale=1]

% \node[inner sep=2pt, circle] (13) at (0,0) [draw] {};
\node[inner sep=2pt, circle] (12) at (-1,1) [draw] {};
\node[inner sep=2pt, circle] (9) at (0,2) [draw] {};
\node[inner sep=2pt, circle] (11) at (-2,2) [draw] {};
\node[inner sep=2pt, circle] (4) at (-1,3) [draw] {};
% \node[inner sep=2pt, circle] (10) at (-3,3) [draw] {};
\node[inner sep=2pt, circle] (8) at (-4,4) [draw] {};
\node[inner sep=2pt, circle] (0) at (-5,5) [draw] {};
% \node[inner sep=2pt, circle] (7) at (-3,5) [] {};
\node[inner sep=2pt, circle] (6) at (-2,6) [draw] {};
\node[inner sep=2pt, circle] (5) at (-3,7) [draw] {};
\node[inner sep=2pt, circle] (3) at (-1,7) [draw] {};
\node[inner sep=2pt, circle] (1) at (0,8) [draw] {};
\node[inner sep=2pt, circle] (2) at (-2,8) [draw] {};

% \path[style=semithick] (12) edge node[anchor=north]{{}}(13);
\path[style=semithick] (12) edge node[anchor=west]{{}}(9);
\path[style=semithick] (12) edge node[anchor=south]{{}}(11);
\path[style=semithick] (11) edge node[anchor=south]{{}}(4);
\path[style=semithick] (11) edge node[anchor=south]{{}}(8);
% \path[style=semithick] (8) edge node[anchor=south]{{}}(10);
\path[style=semithick] (8) edge node[anchor=south]{{}}(0);
\path[style=semithick] (8) edge node[anchor=south]{{}}(6);
% \path[style=semithick] (7) edge node[anchor=south]{{}}(6);
\path[style=semithick] (6) edge node[anchor=south]{{}}(5);
\path[style=semithick] (6) edge node[anchor=south]{{}}(3);
\path[style=semithick] (3) edge node[anchor=south]{{}}(1);
\path[style=semithick] (3) edge node[anchor=south]{{}}(2);

% \node[anchor = north]  at (13) {13L};
\node[anchor = east]  at (12) {12L};
\node[anchor = east]  at (11) {11L};
% \node[anchor = east]  at (10) {10L};
\node[anchor = east]  at (8) {8L};
\node[anchor = east]  at (0) {0L};
\node[anchor = west]  at (4) {4R};
\node[anchor = west]  at (9) {9R};
% \node[anchor = west]  at (7) {7R};
\node[anchor = west]  at (6) {6R};
\node[anchor = west]  at (3) {3R};
\node[anchor = west]  at (1) {1R};
\node[anchor = east]  at (2) {2L};
\node[anchor = east]  at (5) {5L};
\end{tikzpicture}
$$
In the first step, we construct the induced discrete Morse function on a path of the underlying Morse labeled merge tree:
$$
\begin{tikzpicture}
    \node[inner sep=2pt, circle] (0) at (0,0) [draw] {};
    \node[inner sep=2pt, circle] (5) at (2,0) [draw] {};
    \node[inner sep=2pt, circle] (2) at (4,0) [draw] {};
    \node[inner sep=2pt, circle] (1) at (6,0) [draw] {};
    \node[inner sep=2pt, circle] (4) at (8,0) [draw] {};
    \node[inner sep=2pt, circle] (9) at (10,0) [draw] {};

    \node[anchor = north]  at (0) {0};
    \node[anchor = north]  at (5) {5};
    \node[anchor = north]  at (2) {2};
    \node[anchor = north]  at (1) {1};
    \node[anchor = north]  at (4) {4};
    \node[anchor = north]  at (9) {9};

    \path[style=semithick] (0) edge node[anchor=south]{{8}}(5);
    \path[style=semithick] (5) edge node[anchor=south]{{6}}(2);
    \path[style=semithick] (2) edge node[anchor=south]{{3}}(1);
    \path[style=semithick] (1) edge node[anchor=south]{{11}}(4);
    \path[style=semithick] (4) edge node[anchor=south]{{12}}(9);
\end{tikzpicture}
$$
In the second step, we add the closing edges, which correspond to the cycle nodes of the given gML tree.
$$
\begin{tikzpicture}
    \node[inner sep=2pt, circle] (0) at (0,0) [draw] {};
    \node[inner sep=2pt, circle] (5) at (2,0) [draw] {};
    \node[inner sep=2pt, circle] (2) at (4,0) [draw] {};
    \node[inner sep=2pt, circle] (1) at (6,0) [draw] {};
    \node[inner sep=2pt, circle] (4) at (8,0) [draw] {};
    \node[inner sep=2pt, circle] (9) at (10,0) [draw] {};

    \node[anchor = north]  at (0) {0};
    \node[anchor = north]  at (5) {5};
    \node[anchor = north]  at (2) {2};
    \node[anchor = north]  at (1) {1};
    \node[anchor = north]  at (4) {4};
    \node[anchor = north]  at (9) {9};

    \path[style=semithick] (0) edge node[anchor=south]{{8}}(5);
    \path[style=semithick] (5) edge node[anchor=south]{{6}}(2);
    \path[style=semithick] (2) edge node[anchor=south]{{3}}(1);
    \path[style=semithick] (1) edge node[anchor=south]{{11}}(4);
    \path[style=semithick] (4) edge node[anchor=south]{{12}}(9);

    \draw (0) .. controls (1,-1) .. (5) node[midway, below] {10};
    \draw (5) .. controls (3,-1) .. (2) node[midway, below] {7};
    \draw (4) .. controls (9,-1) .. (9) node[midway, below] {1};
\end{tikzpicture}
$$
\end{example}
\begin{remark}
    We compare the discrete Morse fucntion on a graph from \Cref{ex: first dmf} and the induced discrete Morse function, see \Cref{ex:induced dMf}, of its induced generalized Morse labeled merge tree, see \Cref{ex: generalized Morse merge tree}: at first, it might seem that the two graphs look quite different. As it turns out, the two graphs are related to each other by removing and reattaching closing edges to the connected components of sublevel complexes they are connected with. Thus, the induced generalized Morse labeled merge tree only ``sees'' which connected components merge and which connected components contain cycles along the filtration, but not the exact attaching information. We make this relation precise in \Cref{cmequiv} and \Cref{DmfongraphvsMltree}. 
\end{remark}
\begin{lemma}\label{}
We have $M(\Bar{X},f_{\lvert \Bar{X}})\cong \Bar{M}(X,f)$ as Morse labeled merge trees. 
\end{lemma}
\begin{proof}
The construction of the induced generalized merge tree induces a bijection $M\colon X \rightarrow V(M(X,f))$. It follows immediately by construction that $M$ bijectively maps closing edges to nodes of $ C(M(X,f))$. Hence removing the closing edges from $(X,f)$, that is, passing on to $(\Bar{X},f)$, precisely removes the nodes of $C(M(X,f))$, which corresponds to passing on to $\Bar{M}(X,f)$. Hence, the statement holds because the values of $f$ on non-closing edges are not changed.
\end{proof}
\begin{definition}[{\cite[Definition 2.42]{brueggemann2021merge}}]\label{symm}
     Let $f\colon X \rightarrow \mathbb{R}$ be a dMf on a graph. For each non-empty connected component $X_c[v]$ of a sublevel complex $X_c$ we denote by $\operatorname{Aut}(X_c[v])$ the group of simplicial automorphisms of $X_c[v]$. Each $\xi \in \operatorname{Aut}(X_c[v])$ can be extended by the identity to a set function that is a self-bijection $X\rightarrow X$. The group $\widetilde{\operatorname{Aut}}(X_c[v])$ is defined to be the group of said extensions of elements of $\operatorname{Aut}(X_c[v])$ by the identity\footnote{That is, elements of $\operatorname{Aut}(X_c[v])$ are self-bijections of $X$ that restrict to a simplicial automorphism on $X_c[v])$ and to the identity on $X\setminus X_c[v])$.}. The group operation on $\widetilde{\operatorname{Aut}}(X_c[v])$ is the composition of self-bijections of $X$. We call the elements of $\widetilde{\operatorname{Aut}}(X_c[v])$ \textit{elementary sublevel automorphisms}.
     % The total order on $\operatorname{Cr}(f)$ induced by $f$ induces chains  $\widetilde{\operatorname{Aut}}(X^f_{c_0}[v])\subset \widetilde{\operatorname{Aut}}(X^f_{c_1}[v]) \subset \dots$ of inclusions of subgroups. Moreover, we have inclusions $\widetilde{\operatorname{Aut}}(X^f_{c_i}[v])\subset \widetilde{\operatorname{Aut}}(X^f_{c_j}[v])=\widetilde{\operatorname{Aut}}(X^f_{c_j}[w]) \supset \widetilde{\operatorname{Aut}}(X^f_{c_i}[w])$ if $v$ and $w$ are in different connected components of some sublevel complex $X^f_{c_i}$ that merge together in some other sublevel complex $X^f_{c_j}$ for $j>i$.  
     We define the \textit{sublevel automorphism group} of $(X,f)$, denoted by $\operatorname{Aut}_{sl}(X,f)$, as 
     $$\operatorname{Aut}_{sl}(X,f) \coloneqq \freeprod\limits_{c \in \operatorname{Cr}(f),v \in X} \widetilde{\operatorname{Aut}}(X_c[v])/\sim,$$
     where $\ast$ denotes the free product of groups and $\sim$ is defined by
     $$ \xi\xi'b\sim \begin{cases}
         \xi\circ \xi' & \text{if }a,b\in \widetilde{\operatorname{Aut}}(X_c[v]) \text{ for the same }X_c[v]\\
         \xi'\xi & \text{if } \xi,\xi' \text{ belong to different connected components of sublevel complexes}
     \end{cases}$$
     We call the elements of $\operatorname{Aut}_{sl}(X,f)$ \textit{sublevel automorphisms}.  
\end{definition}

Note that an element $\xi\in \operatorname{Aut}_{sl}(X,f)$ is not necessarily an automorphism of $X$, but only a collection of self-bijections of $X$ that restrict to automorphisms on certain connected components of some sublevel complex $X_c[v]$, and the identity outside of $X_c[v]$. 
% It is immediate that the definition of sublevel automorphism groups can be extended to arbitrary cellular filtrations of CW complexes. \par
Furthermore, notice that $\operatorname{Aut}_{sl}(X,f)$ is by definition isomorphic to the product of the $\widetilde{\operatorname{Aut}}(X_c[v])$. We chose to phrase $\operatorname{Aut}_{sl}(X,f)$ as a quotient of a free product to provide more clarity in \Cref{Lem:DMFCompositeSublevelAutoIsDMF} and \Cref{Prop:ActionOfSublevelAutos}. Moreover, $\operatorname{Aut}_{sl}(X,f)$ is not a subgroup of the group of self-bijections of $X$. 

% \begin{lemma}\label{lemma sublevel automorphisms cummutative}
%     Sublevel automorphisms of different levels commute with each other.
% \end{lemma}
% \begin{proof}
%     Let $(X,f)$ be a dMf on a graph and let $a_1,a_2\in \operatorname{Aut}_{sl}(X,f)$ be sublevel automorphisms. If $a_1$ and $a_2$ act on disjoint connected components of sublevel complexes, there is nothing to show. Without loss of generality, assume that $a_2$ acts on $X^f_{c_2}[v]$ and $a_1$ acts on $X^f_{c_1}[v]\subset X^f_{c_2}[v]$, i.e. $a_2$ is of a higher level than $a_1$ and $v$ is a vertex that both $a_1$ and $a_2$ act on. We observe that applying $a_2$ also affects $X^f_{c_1}[v]\subset X^f_{c_2}[v]$ and hence changes the action of $a_1$ on $X$. The other way around, applying $a_1$ changes how $X^f_{c_1}[v]$ is embedded into $X^f_{c_2}[v]$ and hence how $a_2$ acts on $X^f_{c_2}[v]$. Since the way $a_1$ and $a_2$ affect each other is governed by the inclusion $X^f_{c_1}[v]\subset X^f_{c_2}[v]$, the two different sublevel automorphisms change each other exactly as needed to be commutative.
% \end{proof}

\begin{example}\label{examplesofslsymm}
    We consider three instructive examples of sublevel-automorphism groups:
    \begin{enumerate}
        \item The path with $n$ vertices $P$,
        \item the star graph with $n+1$ vertices $S$, and
        \item the cycle graph with $n$ vertices $C$.
        \end{enumerate}
        \textbf{(1)}\par
        Consider the path $P$ with $n$ vertices with a critical dMf defined by counting from left to right:
        \begin{center}
        \begin{tikzpicture}
        \node[inner sep=2pt, circle] (0) at (0,0) [draw] {};
        \node[inner sep=2pt, circle] (1) at (3,0) [draw] {};
        \draw (9.5,0) node {$\dots$};
        \node[inner sep=2pt, circle] (3) at (13,0) [draw] {};
        \node[inner sep=2pt, circle] (4) at (16,0) [draw] {};
        \node[inner sep=2pt, circle] (5) at (6,0) [draw] {};

        \node[anchor = north]  at (0) {{$0$}};
        \node[anchor = north]  at (1) {{$1$}};
        \node[anchor = north]  at (3) {{$2(n-1)-3$}};
        \node[anchor = north]  at (4) {{$2n-3$}};
        \node[anchor = north]  at (5) {{$3$}};
        \draw[-]  (0)--(1) node[midway, above] {2};
        \draw[-]  (5)--(8.5,0);
        \draw[-]  (3)--(4) node[midway, above] {2n-2};
        \draw[-]  (10.5,0)--(3);
        \draw[-]  (1)--(5) node[midway, above] {4};
        \end{tikzpicture}
        \end{center}
        We observe that each connected component of a sublevel set only has either the trivial group or a group generated by exactly one reflection as its automorphism group.
        % The reflections from the different levels commute with one another by \Cref{lemma sublevel automorphisms cummutative}. 
        Hence, we have $\operatorname{Aut}(P^f_{f(e)}[e])\cong \Sigma_2$ for all edges $e$ and $\operatorname{Aut}_{sl}(P,f)\cong \prod\limits_{k=1}^{n-1} \Sigma_2$, where $\Sigma_k$ denotes the symmetric group on $k$ elements. This way, we realized the Young subgroup $\prod\limits_{k=1}^{n-1} \Sigma_2\subset \Sigma_{2(n-1)}$ as a group of sublevel automorphisms of a filtered space with the associated constant sequence of automorphism groups of sublevel complexes $\operatorname{Aut}(P^f_{f(e)}[e])\cong \Sigma_2$.\\
        \textbf{(2)}\par
        We consider the star graph with $n+1$ vertices $S$ together with a critical dMf that attains its minimum at the center and otherwise assigns values pairwise to the outer vertices and their adjacent edges:
        \begin{center}
            \begin{tikzpicture}
                \node[inner sep=2pt, circle] (0) at (0,0) [draw] {};
                \node[inner sep=2pt, circle] (1) at (1.41,1.41) [draw] {};
                \node[inner sep=2pt, circle] (2) at (2,0) [draw] {};
                \node[inner sep=2pt, circle] (3) at (1.41,-1.41) [draw] {};
                \draw (-1.41,-1.41) node {$\dots$};
                \node[inner sep=2pt, circle] (4) at (-1.41,1.41) [draw] {};

                \node[anchor = east]  at (0) {{$0$}};
                \node[anchor = south]  at (1) {{$1$}};
                \node[anchor = west]  at (2) {{$3$}};
                \node[anchor = north]  at (3) {{$5$}};
                \node[anchor = south]  at (4) {{$n$}};

                \draw[-]  (0)--(1) node[midway, above] {2};
                \draw[-]  (0)--(2) node[midway, above] {4};
                \draw[-]  (0)--(3) node[midway, below] {6};
                \draw[-]  (0)--(4) node[midway, left] {$n+1$};
            \end{tikzpicture}
        \end{center}
        Let $v_k$ be the vertex with label $k$. It is immediate that $S_{2k+2}^f[v_{2k+1}]$ is the star graph with $k+2$ vertices. Hence, we have $Aut(S_{2}^f[v_{1}])\cong \Sigma_{2}$ and $\mathrm{Aut}(S_{2k+2}^f[v_{2k+1}])\cong \Sigma_{k+1}$ for $k\geq 1$, where $\Sigma_k$ denotes the symmetric group on $k$ elements. We have $\mathrm{Aut}_{sl}(S,f)\cong \Sigma_2 \times \prod\limits_{k=2}^n \Sigma_k$. So we realized $\Sigma_2 \times \prod\limits_{k=2}^n \Sigma_k$ as a group of sublevel automorphisms of a filtered space with the associated sequence of groups of automorphisms of connected components of sublevel complexes $\mathrm{Aut}(S_{2}^f[v_{1}])\cong \Sigma_{2}$ and $\mathrm{Aut}(S_{2k+2}^f[v_{2k+1}])\cong \Sigma_{k+1}$ for $k\geq 1$.\\
        \textbf{(3)}\par
        Consider the cycle graph with $n$ vertices $C$ together with the following critical dMf:
        \begin{center}
            \begin{tikzpicture}
                \node[inner sep=2pt, circle] (0) at (0,2) [draw] {};
                \node[inner sep=2pt, circle] (1) at (1.41,1.41) [draw] {};
                \node[inner sep=2pt, circle] (2) at (2,0) [draw] {};
                \node[inner sep=2pt, circle] (3) at (1.41,-1.41) [draw] {};
                \node[inner sep=2pt, circle] (4) at (0,-2) [draw] {};
                \node[inner sep=2pt, circle] (5) at (-1.41,-1.41) [draw] {};
                \node[inner sep=2pt, circle] (6) at (-2,0) [draw] {};
                \node[inner sep=2pt, circle] (7) at (-1.41,1.41) [draw] {};

                \node[anchor = south]  at (0) {{$0$}};
                \node[anchor = west]  at (1) {{$1$}};
                \node[anchor = west]  at (2) {{$3$}};
                \node[anchor = west]  at (3) {{$5$}};
                \node[anchor = north]  at (4) {{$7$}};
                \node[anchor = east]  at (5) {{$2n-7$}};
                \node[anchor = east]  at (6) {{$2n-5$}};
                \node[anchor = east]  at (7) {{$2n-3$}};

                \draw[-]  (0)--(1) node[midway, above] {$2$};
                \draw[-]  (1)--(2) node[midway, right] {$4$};
                \draw[-]  (2)--(3) node[midway, right] {$6$};
                \draw[-]  (3)--(4) node[midway, below] {$8$};
                \draw[dotted]  (4)--(5);
                \draw[-]  (5)--(6) node[midway, left] {$2n-4$};
                \draw[-]  (6)--(7) node[midway, left] {$2n-2$};
                \draw[-]  (7)--(0) node[midway, above] {$2n$};
            \end{tikzpicture}
        \end{center}
        Let $e$ be any edge except for the one labeled $2n$. Then we have $\operatorname{Aut}(P^f_{f(e)}[e])\cong \Sigma_2$ as in the first example. Let $\Tilde{e}$ be the edge labeled $2n$. Then we have $\operatorname{Aut}(P^f_{f(\Tilde{e})}[\Tilde{e}])\cong D_{n}$ where $D_{n}$ denotes the dihedral group of order $2n$, i.e. the symmetries of the regular $n$-gon. This leads us to $\mathrm{Aut}_{sl}(S,f)\cong D_n \times \prod\limits_{k=1}^{n-1} \Sigma_2$. This way, we realized $D_n \times \prod\limits_{k=1}^{n-1} \Sigma_2$ as a group of sublevel automorphisms of a filtered space with the associated constant sequence of automorphism groups of sublevel complexes $\operatorname{Aut}(C^f_{f(e)}[e])\cong \Sigma_2$ for $f(e)\leq 2n-2$ and $\operatorname{Aut}(C^f_{2n}[\Tilde{e}])\cong D_n$.
\end{example}

\begin{definition}\label{Def:SymmEquivalence}
     Let $f\colon X\rightarrow \mathbb{R}$ and $g\colon X\rightarrow \mathbb{R}$ be dMfs on a graph $X$. We call $f$ and $g$ \textit{sublevel-equivalent} if they have the same critical values and isomorphic sublevel complexes. If additionally $g=f\circ \xi$ holds for a sublevel automorphism $\xi \in \operatorname{Aut}_{sl}(X,f)$, then we call $f$ and $g$ \textit{symmetry-equivalent}. We call the map $\xi$ a \textit{symmetry equivalence} from $f$ to $g$. \par
     We call two dMfs $f\colon X\rightarrow \mathbb{R}$ and $g\colon Y\rightarrow \mathbb{R}$ \textit{symmetry-equivalent} if there is a simplicial isomorphism $\varphi\colon X\rightarrow Y$ such that $f$ and $g\circ \varphi$ are symmetry-equivalent.
\end{definition}
Having these definitions established, we are able to consider the action of $\operatorname{Aut}_{sl}(X,f)$ on the symmetry equivalence class of $(X,f)$. 
% In particular, we can consider the partial group action of the union of all sublevel automorphism groups, as a subgroup of the group of self-bijections on $X$, on the space of dMfs on $X$.
\begin{remark}\label{Rem:sublevelAutomorphismsCommute}
    If two dMfs on graphs  $f\colon X\rightarrow \mathbb{R}$ and $g\colon Y\rightarrow \mathbb{R}$ are symmetry-equivalent, then their sublevel automorphism groups are isomorphic because the two dMfs induce isomorphic filtrations. \par
    We want to remark at this point that even though the elements of sublevel automorphism groups $\operatorname{Aut}_{sl}(X,f)$ are collections of self-bijections of the respective graph $X$, the group structure is different from the group of self-bijections of $X$. In particular, the group structure of $\operatorname{Aut}_{sl}(X,f)$ is constructed such that symmetry equivalences of different levels always commute with each other although the corresponding self-bijections of $X$ do not necessarily commute. This is necessary for the desired action on the set of dMfs on $X$: symmetry equivalences exist due to the existence of different filtrations with isomorphic sublevel complexes that are embedded differently into $X$. That is, if $\varphi$ and $\psi$ are symmetry equivalences that belong to two connected components of sublevel complexes such that one is contained in the other. Without loss of generality, $\varphi$ is of a higher level than $\psi$, then $\varphi$ changes the location of the connected component that belongs to $\psi$ precisely such that the actions of $\varphi$ and $\psi$ on the set of dMfs commute.
\end{remark}
\begin{lemma}\label{Lem:DMFCompositeSublevelAutoIsDMF}
    Let $f\colon X \rightarrow \mathbb{R}$ be a dMf on a graph $X$ and let $\xi\in \widetilde{\operatorname{Aut}}(X_c[v]) \subset  \operatorname{Aut}_{sl}(X,f)$ be an elementary sublevel automorphism. Then $f\circ \xi$ is a dMf on $X$, which is symmetry equivalent to $f$.\par  Moreover, $\xi$ induces an isomorphism $\operatorname{Aut}_{sl}(X,f)\cong\operatorname{Aut}_{sl}(X,g)$ by precomposition. 
\end{lemma}
\begin{proof}
    In order to prove that function $f \circ \xi$ is a dMf, we note that $\xi$ is a simplicial automorphism on $X_c[v]$ and the identity outside of $X_c[v]$. Since $\xi$ is in particular a self-bijection of $X$, $f\circ \xi$ is still at most 2-1. Due to $X_c[v]$ being contained in a sublevel complex, all values of $f$ outside of $X_c[v]$ are strictly larger than the ones inside $X_c[v]$. In particular, $f$ is strictly monotone on all face relations at the boundary of $X_c[v]$, i.e., between simplices of $X_c[v]$ and simplices outside $X_c[v]$. Thus, the action of $\xi$ outside of $X_c[v]$ does not affect monotonicity and generacy. Inside $X_c[v]$, $\xi$ acts as a simplicial automorphism, hence $f\circ \xi$ also satisfies monotonicity and generacy. Furthermore, $f\circ \xi$ is symmetry-equivalent to $f$ by \Cref{Def:SymmEquivalence}.\par
    For the second statement, note that $\xi$ induces an isomorphism between the filtrations induced by $f,f\circ \xi$, respectively. That is, $\xi$ bijectively maps connected components of sublevel complexes of $f$ to connected components of sublevel complexes of $g$ in an inclusion and filtration preserving way. Hence, it follows directly from the presentation of $\operatorname{Aut}_{sl}$ given in \Cref{symm} that $\xi$ is an isomorphism.
\end{proof}
\begin{proposition}\label{Prop:ActionOfSublevelAutos}
    Let $X$ be a graph, let $f\colon X \rightarrow \mathbb{R}$ be a discrete Morse function, and let $\operatorname{Aut}_{sl}(X,f)$ be the group of sublevel automorphisms of $(X,f)$. Then there is an action of $\operatorname{Aut}_{sl}(X,f)$ on the symmetry equivalence class of $f$ given by: for any dMf $g$ that is symmetry-equivalent to $f$, and any elementary sublevel equivalence $\xi$, we define $g*\xi\coloneqq g\circ \tilde{\xi}$, where $\tilde{\xi}\in \operatorname{Aut}_{sl}(X,g)$ is the elementary sublevel equivalence that corresponds to $\xi$ under the automorphism $\operatorname{Aut}_{sl}(X,f)\cong\operatorname{Aut}_{sl}(X,g)$ from \Cref{Lem:DMFCompositeSublevelAutoIsDMF}. For arbitrary elements of $\operatorname{Aut}_{sl}(X,f)$, the group action is defined by the successive action of elementary sublevel automorphisms.
\end{proposition}
\begin{proof}
    It follows from successive application of \Cref{Lem:DMFCompositeSublevelAutoIsDMF} that $g*\xi$ is well defined for any dMf $g$ that is symmetry-equivalent to $f$, and any elementary equivalence $\xi$. The compatibility, i.e., that for any $\xi,\xi'\in\operatorname{Aut}_{sl}(X,f)$, we have $g*\xi*\xi'=g*(\xi\cdot \xi')$, where $\cdot$ denotes the multiplication in $\operatorname{Aut}_{sl}(X,f)$ follows by construction of the group action and $\operatorname{Aut}_{sl}(X,f)$, and by \Cref{Rem:sublevelAutomorphismsCommute}.
    % Without loss of generality, we can show this for the case of elementary sublevel automorphisms because the general statement then follows by successive appplication of the case of elementary sublevel automorphisms.\par
    % If $\xi$ and $\xi'$ are both elementary sublevel equivalences of the same connected component of a sublevel complex 
\end{proof}
Next, we introduce the more general notion of component-merge equivalence.
\begin{definition}[{\cite[Definition 2.50]{brueggemann2021merge}}]
\label{cmequiv}
	Let $(X,f)$ and $(X',f')$ be critical dMfs on connected graphs. A \textit{component-merge equivalence (CM equivalence)} of level $a$\footnote{We emphasize the notion of the level of CM equivalences in order to highlight the recursive nature of this definition. In situations when the specific level of a CM equivalence is not of importance, we sometimes drop the level in the notation.} is a bijection $\varphi \colon X \rightarrow X'$ such that at least one of the following two cases holds:
	\begin{enumerate}
		\item $\varphi$ is a symmetry equivalence that involves sublevel automorphisms of at most level $a$.
		\item $\varphi$ fulfills the following three conditions:
	\begin{itemize}
		\item $f' \circ \varphi=f$,
		\item $\varphi$ induces a bijection between the sets of connected components of sublevel complexes such that the restriction $\varphi_{\lvert X_{a-\varepsilon}[v]}\colon X_{a-\varepsilon}[v] \rightarrow X'_{a-\varepsilon}[\varphi(v)]$ to each connected component is a CM equivalence of a level $\leq a$, and
		\item the edge $\sigma \in X$ with $f(\sigma)=a$ merges two connected components $X_{a-\varepsilon}[v_1]$ and $X_{a-\varepsilon}[v_2]$ in $X_{a}[v_1]=X_{a}[v_2]$ if and only if the edge $\varphi(\sigma)$ merges the corresponding two connected components $X'_{a-\varepsilon}[\varphi(v_1)]$ and $X'_{a-\varepsilon}[\varphi(v_2)]$ in $X'_{a}[\varphi(v_1)]=X'_{a}[\varphi(v_2)]$. Otherwise, if the edge $\sigma \in X$ with $f(\sigma)=a$ does not merge two connected components but rather closes a circle within a connected component $X_{a-\varepsilon}[v]$, then and only then $\varphi(\sigma)$ closes a circle within $X'_{a-\varepsilon}[\varphi(v)]$. 
	\end{itemize}
If $\varphi$ re-attaches the critical edge labeled $a$,
we call $\varphi$ \textit{non-trivial}. Moreover, if $\varphi$ re-attaches the critical edge of level $a$ and acts as a symmetry equivalence everywhere else, we say that $\varphi$ is \textit{elementary of level $a$}. If $\varphi$ does not re-attach any critical edge, i.e., if $\varphi$ is a symmetry equivalence, we call $\varphi$ a \textit{trivial CM equivalence}.
\end{enumerate}

\begin{remark}\label{extend}
Extending the notion of CM equivalences to dMfs with matched cells is a bit tedious. We would like to suggest getting rid of matched cells by identifying arbitrary dMfs on graphs with critical dMfs on the corresponding graph that arises by collapsing matched cells beforehand even though self-loops might arise in this process. Nonetheless, the newly created self-loops are critical by construction and the definition above works in this context. 
\end{remark}

\begin{example}
Let $f\colon X \to \mathbb{R}$ be the complex with dMf on the left and $f'\colon X'\to \mathbb{R}$ be the complex with dMf on the right.
$$
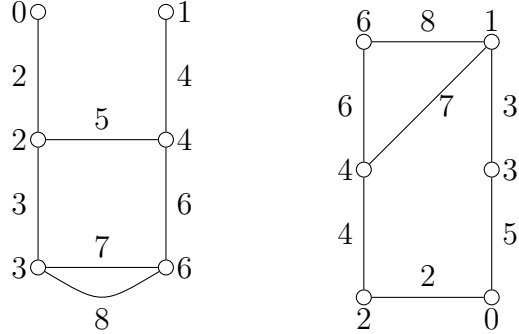
\begin{tikzpicture}[scale=1.7]

\node[inner sep=2pt, circle] (3) at (0,0) [draw] {};
\node[inner sep=2pt, circle] (2) at (0,1) [draw] {};
\node[inner sep=2pt, circle] (0) at (0,2) [draw] {};
\node[inner sep=2pt, circle] (6) at (1,0) [draw] {};
\node[inner sep=2pt, circle] (4) at (1,1) [draw] {};
\node[inner sep=2pt, circle] (1) at (1,2) [draw] {};

\draw[-]  (3)--(6) node[midway, above] {$7$};
\draw[-]  (3)--(2) node[midway, left] {$3$};
\draw[-]  (2)--(0) node[midway, left] {$2$};
\draw[-]  (6)--(4) node[midway, right] {$6$};
\draw[-]  (4)--(1) node[midway, right] {$4$};
\draw[-]  (2)--(4) node[midway, above] {$5$};
\draw (3) .. controls (.5,-.3)  .. (6) node[midway, below] {$8$};

\node[anchor = east]  at (3) {{$3$}};
\node[anchor = west]  at (6) {{$6$}};
\node[anchor = east]  at (2) {{$2$}};
\node[anchor = west]  at (4) {{$4$}};
\node[anchor = east]  at (0) {{$0$}};
\node[anchor = west]  at (1) {{$1$}};

\end{tikzpicture}
\qquad
\qquad
\begin{tikzpicture}[scale=1.7]

\node[inner sep=2pt, circle] (2) at (0,0) [draw] {};
\node[inner sep=2pt, circle] (4) at (0,1) [draw] {};
\node[inner sep=2pt, circle] (6) at (0,2) [draw] {};
\node[inner sep=2pt, circle] (0) at (1,0) [draw] {};
\node[inner sep=2pt, circle] (3) at (1,1) [draw] {};
\node[inner sep=2pt, circle] (1) at (1,2) [draw] {};

\draw[-]  (2)--(0) node[midway, above] {$2$};
\draw[-]  (2)--(4) node[midway, left] {$4$};
\draw[-]  (4)--(6) node[midway, left] {$6$};
\draw (3) .. controls (.75,1.5)  .. (1) node[midway, left] {$7$};
\draw[-]  (0)--(3) node[midway, right] {$5$};
\draw[-]  (6)--(1) node[midway, above] {$8$};
\draw[-]  (3)--(1) node[midway, right] {$3$};

\node[anchor = west]  at (3) {{$3$}};
\node[anchor = south]  at (6) {{$6$}};
\node[anchor = north]  at (2) {{$2$}};
\node[anchor = east]  at (4) {{$4$}};
\node[anchor = north]  at (0) {{$0$}};
\node[anchor = south]  at (1) {{$1$}};

\end{tikzpicture}
$$
Then a CM equivalence $\varphi\colon X\to X'$ of critical levels $a=7$ is given by $\varphi(v)=v'$ whenever $f(v)=f'(v')$ on vertices and $\varphi(e)=e'$ whenever $f(e)=f'(e')$ on edges. We remark that according to \Cref{extend} matched simplices can be arbitrarily added and removed from connected components by CM equivalences. After performing the prescribed collapse, it becomes clear that $\varphi$ only re-attaches the edge labeled $7$ from an edge between the critical vertices labeled $0$ and $1$ to a critical self-loop at the critical vertex labeled $1$.
\end{example}

\end{definition}
\begin{remark}\label{cmclo}
It is clear from the case distinction made in  \Cref{cmequiv} that any CM equivalence $\varphi \colon (X,f) \rightarrow (X',f')$ restricts to a bijection $\varphi_{\lvert C(X,f)}\colon C(X,f)\rightarrow C(X',f')$.
\end{remark}

At the end of this section, we want to provide a different point of view on CM equivalences. As opposed to the case of symmetry equivalences, we cannot describe the action of the group of CM equivalences on the CM equivalence class of some dMf on a graph in terms of some group action on a space because CM equivalences change the space at hand. \par
Instead we propose to consider this operation as a digraph which has the CM equivalence class of a dMf on a graph $(X,f)$ as vertices and elementary CM equivalences, i.e., ones that are either a symmetry equivalence of only one connected component or a CM equivalence of some level $a$, as edges. In the previous proof, we already used the fact that CM equivalences can be decomposed into a sequence of CM equivalences of separate levels $a$. Such CM equivalences of level $a$ are determined by which edge $e$ they reattach, that is, which vertices the boundary vertices of $e$ are swapped with. This allows us to order the outgoing edges at each vertex linearly.\par
We identify edges $e$ with the ordered pair of the Morse labels of their boundary vertices $(c_1,c_2)$ with the convention that the smaller label always comes first, that is, $c_1<c_2$. Let $\varphi$ be a CM equivalence of level $f(e)$ that maps $(c_1,c_2)\mapsto (c_1',c_2')$. Then we label $\varphi$ with the ordered tuple $(c_1,c_2,c_1',c_2')$ and order the outgoing CM equivalences by level and among the same level by the lexicographical ordering on these labels. We define trivial CM equivalences, i.e., symmetry equivalences, to be less then non-trivial ones and order them by level and inside their level by minimal label on the connected component involved. For the same connected component at the same level, we define the order given by a lexicographic order on words which describe the symmetries similar to the case of CM equivalences. Thus, we can explore the CM equivalence class of any dMf on a graph with any standard exploration algorithm for any edge labeled digraph. 
\begin{remark}
    The mentioned point of view on CM equivalence classes can be phrased in a more categorical language: the described digraph encodes the data of a groupoid that describes the action of CM equivalences on their corresponding CM equivalence class. 
\end{remark}
\section{Inverse Problem for Multigraphs}
In this section we want to describe the relationship between dMfs on graphs, gML trees, and generalized merge trees. The results are summarized in \Cref{diagramminverseproblem}.\\
\begin{figure}[H]
\centering
\begin{tikzpicture}[scale=1]
\node (Mer) at (0,4-.3) {$gMer$};
\node (DMF) at (8,4-.3) {$DMF^{\text{crit}}_{graphs}$};
% \node (MoT) at (-0.2,2+.3) {$gMOT$};
\node (MlT) at (4,2+.3) {$gMLT$};
% \draw[->] (0+.6,2+.1+.3) -- (8-.6,2+.1+.3) node[midway,above] {$iML$};
% \draw[->] (8-.6,2-.1+.3) -- (0+.6,2-.1+.3) node[midway,below] {$iMO$};
\draw[->] (8-1,4-.4) -- (0+0.6,4-.4) node[midway, below] {$M(\_ \ ,\_)$};
\draw[->] (0+0.6,4-0.2) -- (8-1,4-0.2) node[midway, above] {$\Phi \circ iML \circ \leq_{sc}$};
\draw[->] (4.8,2.3) -- (8,3.3) node[midway,below] {$\Phi$};
\draw[->] (8-.1,3.4) -- (4.7,2.4) node[midway,above] {$M(\_ \ , \_) \hspace{20pt}$};
\draw[->] (3.3,2.4) -- (0.1,3.5) node[midway,above] {$\operatorname{forget}$};
\draw[->] (0,3.4) -- (3.2,2.3) node[midway,below] {$\lambda_{sc}$};
\end{tikzpicture}
\caption{Relationships between dMfs and merge trees}
    \label{diagramminverseproblem}
\end{figure}
Here, $gMer$ denotes the set of isomorphism classes of generalized merge trees.

\begin{theorem}\label{DmfongraphvsMltree}
Let $DMF^{crit}_{graphs}$ denote the set of CM equivalence classes of dMfs with only critical cells on multigraphs. Let $gMLT$ denote the set of isomorphism classes of gML trees. Then the induced dMf $\Phi$, \Cref{idmf}, and the induced Morse labeled merge tree $M(\_ \ , \_ )$, \Cref{iMLT}, define maps $M(\_ \ ,  \_ ) \colon DMF^{crit}_{graphs} \leftrightarrow gMLT \colon \Phi$ that are inverse of each other in the sense that:
\begin{enumerate}
    \item for any dMf $(X,f)$ with only critical cells, the dMf $\Phi(M(X,f),\lambda_f)$ is CM-equivalent to $(X,f)$, and 
    \item for any gML tree $(T,\lambda)$, we have $M(\Phi T,f_\lambda)\cong (T,\lambda)$.
\end{enumerate}
\end{theorem}
Before we prove this theorem, we need the following lemma:
\begin{lemma}\label{lemma: CmEquivImpliesIsoIndMerTree}
Let $f\colon X \rightarrow \mathbb{R}$ and $f'\colon X' \rightarrow \mathbb{R}$ be CM-equivalent dMfs on multigraphs. Then $M(X,f)\cong M(X',f')$ holds as gML trees.
\end{lemma}
\begin{proof}
Let $\varphi$ be a CM equivalence $\varphi \colon (X,f) \rightarrow (X',f')$. Since we work with a generic version of dMfs which are at most 2-1, at most one non-trivial elementary CM equivalence of level $a$ can occur for any level $a$ because there is at most one critical edge labeled $a$ in $(X,f),(X',f')$, respectively. Thus, we can decompose any CM equivalence into a sequence $(\varphi_a)_a$ of non-trivial elementary CM equivalences of decreasing levels such that each $\varphi_a$ only changes the attachment of the single edge $\sigma$ with $f(\sigma)=a$ and acts as a symmetry equivalence on the rest of graph and dMf. It suffices to consider a single level $a$ because the statement then follows by induction from highest to lowest over all levels $a$. \par
	For such a non-trivial elementary CM equivalence $\varphi_a$ we consider the step of the construction of the induced Ml trees that considers the critical edge $\sigma$ with $f(\sigma)=a$ and the critical edge $\varphi(\sigma)$. If $\sigma$ is not closing, neither is $\varphi(\sigma)$ by \Cref{cmclo} and the inductive step follows by \cite[Proposition 2.52]{brueggemann2021merge}. In the case that $\sigma$ is closing, so is $\varphi(\sigma)$ and we inductively assume that $\varphi$ induces an isomorphism of induced gML trees everywhere outside the subtree corresponding to the connected component of $X^f_{a-\varepsilon}$ that the edge $\sigma$ with $f(\sigma)=a$ is attached to. That is, on the rest of $M(X,f)$ the map $M(\varphi)$ is a bijection compatible with the chiral child relation onto $M(X',f')$ except possibly for the subtree of $M(X',f')$ which corresponds to the connected component of $X'^{f'}_{a-\varepsilon}$ that the edge $\varphi(\sigma)$ is attached to. 
	
	Since the map $\varphi$ is compatible with the dMfs and because it restricts to a CM equivalence $X^f_{a-\varepsilon}\rightarrow X'^{f'}_{a-\varepsilon}$, the dMf $f$ attains the same minima and maxima on the two relevant connected component of $X^{f}_{a-\varepsilon}$ as $f'$ does on its counterpart of $X'^{f'}_{a-\varepsilon}$ via $\varphi$. Since \Cref{iMLT} only considers which connected component the considered edge is attached to, it makes no difference for the isomorphism type of the induced Ml trees that in general $\sigma$ is attached to said connected component of $X^f_{a-\varepsilon}$ at vertices that do not correspond via $\varphi$ to the ones adjacent to $\varphi(\sigma)$ in $X'^{f'}_{a-\varepsilon}$. Thus, the construction of the induced gML tree produces nodes with the same chirality and label for both induced Ml trees in the steps that consider $\sigma,\varphi(\sigma)$, respectively. By assumption, the restriction $\varphi_{X^{f}_{a-\varepsilon}}\colon X^f_{a-\varepsilon}\rightarrow X'^{f'}_{a-\varepsilon}$ is a symmetry equivalence, so the isomorphism of Ml trees extends to the subtrees that correspond to the respective connected components. 
\end{proof}

\begin{proof}[Proof of \Cref{DmfongraphvsMltree}]
It follows from \Cref{lemma: CmEquivImpliesIsoIndMerTree} that $M(\_ \ ,  \_ )$ is well-defined.
\begin{enumerate}
    \item Let $(X,f)$ be a dMf with only critical cells on a graph $X$. We construct a CM equivalence $\varphi (X,f)\rightarrow \Phi (M(X,f))$ as follows: First we consider the spanning trees induced by $(X,f)$ and $(\Phi (M(X,f)),f_{\lambda_f})$ and show that they are CM equivalent. Then we define $\varphi$ on the closing edges and prove that $\varphi$ is a CM equivalence.\par
    By application of \cite[Theorem 5.6]{brueggemann2021merge} we have a CM equivalence $\tilde{\varphi}\colon (\Bar{X},\Bar{f})\rightarrow (\Phi (M(\Bar{X},\Bar{f})),\Bar{f}_{\lambda_{\Bar{f}}} )$. We extend $\tilde{\varphi}$ to a CM equivalence $\varphi \colon (X,f)\rightarrow \Phi (M(X,f))$ by mapping each closing edge $\sigma \in X$ such that $f(\sigma) =a$ to the unique edge $\sigma'\in\Phi (M(X,f))$ with $f_{\lambda_f}(\sigma')=a$. The edge $\sigma'\in\Phi (M(X,f))$ is closing because $a$ does not appear as a label on $(M(\Bar{X},\Bar{f})),\lambda_{\Bar{f}})\cong (\Bar{M}(X,f),\Bar{\lambda}_f)$ since $a$ is the value of the closing critical edge $\sigma \in X$. Furthermore, the connected component of $X_{a-\varepsilon}$ that $\sigma$ is attached to corresponds to the subtree of $M(X,f)$ that consists of all descendants of $M(\sigma)$. By \Cref{idmf}, the edge $\sigma'$ is attached to the connected component of $\Phi (M(X,f)_{a-\varepsilon}$ that corresponds to said subtree. It follows that $\varphi$ is a CM equivalence.
    \item Let $(T,\lambda)$ be a gML tree. Let $c_0<c_1<\dots <c_n$ be the critical values of $f_{\lambda}$ and let $\sigma_i \in \Phi T$ such that $f_\lambda (\sigma_i)=c_i$. We recall that the induced merge tree $M$ defines in particular a bijection between the critical cells of $\Phi T$ and the nodes of $M(\Phi T,f_\lambda)$. For any cell $\sigma \in \Phi T$, we recall that we denote the node of $M(\Phi T,f_\lambda)$ that corresponds to $\sigma$ by $M(\sigma)$. We also recall that $\Phi$, as constructed in \cite[Definition 3.21]{brueggemann2021merge}, comes with a bijection that we extended to cycle nodes in \Cref{idmf} $\phi\colon V(T) \rightarrow \Phi T$. An isomorphism $(\varphi,\operatorname{id_\mathbb{R}})\colon(T,\lambda)\rightarrow M(\Phi T,f_\lambda)$ is given by $\varphi \coloneqq M\circ \phi^{-1}$. It is immediate that $\varphi$ is a bijection because $M$ and $\phi$ are. Furthermore, $\varphi$ is by construction compatible with the respective Morse labelings. It is only left to show that $\varphi$ is compatible with the chiral child relation and the respective roots.\par
    Consider $\sigma_n\in \Phi T$. For both trees, the cell $\sigma_n$ corresponds to the root of the respective tree. In $M(\Phi T,f_\lambda)$ this is the case because $f_\lambda$ attains its maximum on $\sigma_n$. In $(T,\lambda)$ this holds because $\phi(\sigma_n)$ holds the maximal Morse label $\lambda(\phi(\sigma_n))=c_n$. Thus, the map $\varphi$ maps the root of $(T,\lambda)$ to the root of $M(\Phi T,f_\lambda)$. \par
    For each critical edge $\sigma_i \in \Phi T$ we have one of the two cases:
    \begin{enumerate}
        \item[a)] $\sigma_i$ is closing, or
        \item[b)]$\sigma_i$ is not closing.
    \end{enumerate}
    For case b), the proof is identical to the proof of case (2) of \cite[Theorem 5.4]{brueggemann2021merge}. For case a), let $\sigma_i$ be a closing critical edge. 
    In this case, the compatibility with the chiral child relation follows directly by case 1 of \Cref{iMLT} and the property that only children of generalized merge trees need to have the same chirality as their parent node.
    
\end{enumerate}
\end{proof}
\begin{corollary}\label{dmfvsMlorder}
Since the bijection from \Cref{DmfongraphvsMltree} is compatible with the Morse labels, it induces a bijection $M(\_, \ \_ ) \colon DMF^{crit}_{graphs}/_{\leq} \leftrightarrow gMLT/_{\leq} \colon \Phi$ where $/_{\leq}$ denotes dividing by order equivalence.
\end{corollary}
\begin{definition}
Let $\lambda$ and $\lambda'$ be two Morse labelings on a generalized merge tree $T$. 
A \emph{merge equivalence} $(T,\lambda)\rightarrow (T,\lambda')$ of gML trees is a self-bijection $\psi\colon V(T)\cong V(T)$ such that 
\begin{enumerate}
    \item for each inner node $a$ of $T$, the node $a$ is the maximum of a subtree $T'$ of $T$ with respect to $\lambda$ if and only if $\psi(a)$ is the maximum of $T'$ with respect to $\lambda'$, and
    \item for each leaf $a$ of $T$, the node $a$ is the minimum of a subtree $T'$ of $T$ with respect to $\lambda$ if and only if $\psi(a)$ is the minimum of $T'$ with respect to $\lambda'$.
\end{enumerate}
 We call $\lambda$ and $\lambda'$ \textit{merge equivalent} if there exists a merge equivalence $(T,\lambda)\rightarrow (T,\lambda')$. 
 A merge equivalence $(T,\lambda)\rightarrow (T',\lambda')$ between different gML trees is a concatenation of an isomorphism $\varphi\colon T\rightarrow T'$ of underlying generalized merge trees and a merge equivalence $(T,\lambda)\xrightarrow{\psi}(T,\varphi^*\lambda')\xrightarrow{\varphi}(T',\lambda')$. 
    \end{definition}
    \begin{remark}
        It is straightforward to see that order equivalences are special cases of merge equivalences.
    \end{remark}
\begin{proposition}
Any two Morse labelings $\lambda$ and $\lambda'$ on a generalized merge tree $T$ are merge equivalent.
\end{proposition}
\begin{proof}The statement is proved inductively. Let $a$ be the minimal leaf of a subtree $T'$ of $T$ with respect to $\lambda$. Then $a$ needs to be the minimal leaf of $T'$ with respect to $\lambda'$ because otherwise $\lambda'$ would fail to be a Morse labeling due to \Cref{morseord}. The statement for inner nodes follows similarly. 
\end{proof}
\begin{corollary}\label{corollary merge equivalence}
Two gML trees have isomorphic underlying generalized merge trees if and only if they are merge equivalent. In particular, two (not generalized) ML trees have isomorphic underlying (not generalized) merge trees if and only if they are merge equivalent.
\end{corollary}
For any generalized merge tree $T$, there are several ways to induce canonical\footnote{Here we use the term canonical to refer to Morse labelings that are explicitly constructed from the combinatorial data of the given generalized merge tree, i.e. the parent child relation and the chirality. We use this in contrast to other Morse labelings that are not so easily found.} Morse labelings on $T$. We introduce the sublevel-connected Morse labeling (generalization of \cite[Definition 4.5]{brueggemann2021merge}) on any given generalized merge tree in the following:\par
To define the sublevel-connected Morse labeling, we first observe that every node $a$ of $T$ is uniquely determined by the shortest path from the root to $a$. We recall that the depth of $T$ is the maximal length of any path in $T$ that appears as the shortest path from the root to a leaf. Because $T$ is chiral, we can identify such shortest paths with certain words:
\begin{definition}[{\cite[Definition 3.1]{brueggemann2021merge}}]\label{pw}
	Let $T$ be a generalized merge tree of depth $n$ and let $a$ be a node of $T$. The \textit{path word} corresponding to $a$ is a word $a_0a_1\dots a_n \in \{L,R,\_\}^{n+1}$ where $\_$ denotes the empty letter. If $a$ is of depth $k$, the letters $a_0\dots a_k$ are given by the chirality of the nodes belonging to the shortest directed path from the root to $a$. The letters $a_{k+1}\dots a_n$ are then empty.
\end{definition}
\begin{remark}\label{maxk}
    Let $a,b$ be nodes of a generalized merge tree $T$ and let $a_0a_1\dots a_n$ be the path word corresponding to $a$ and $b_0b_1\dots b_n$ be the path word corresponding to $b$. Then the equation $a_0=b_0=L$ always holds because we consider paths that begin at the root. Because $a_0=b_0=L$ and because we consider finite trees, there is always a maximal $k \in \mathbb{N}$ such that $a_i=b_i$ holds for all $i\leq k$. Furthermore, the last non-empty letter of a path word is always the chirality of the considered node. We have examples of path words in \Cref{ex: sc order example}.
\end{remark}

\begin{definition}[{\cite[Definition 4.1, Definition 4.5]{brueggemann2021merge}}]\label{scmo}
Let $T$ be a generalized merge tree. 
We define the \textit{sublevel-connected Morse labeling $\lambda_{sc}$} on the nodes of $T$ as follows:\par
Let $a,b$ be arbitrary nodes of $T$.
Let $a_0a_1\dots a_n$ be the path word corresponding to $a$ and $b_0b_1\dots b_n$ the path word corresponding to $b$ (see \Cref{pw}).
Furthermore, let $k\in\mathbb{N}$ be maximal such that $a_i=b_i$ for all $i\leq k$. We define $a\leq_{sc} b$ if and only if one of the following cases holds:
\begin{enumerate}[a)]
	\item[a$_1$)] $a_k=b_k=L$, $a_{k+1}=L \text{ and } b_{k+1}=R$ 
    \item[a$_2$)] $a_k=b_k=R$, $a_{k+1}=R \text{ and } b_{k+1}=L$
	\item[b)] $b_{k+1}=\_$
	\item[c)] $a=b$
\end{enumerate}
We define the \emph{sublevel-connected Morse labeling} $\lambda_{sc}$ by mapping the nodes of $T$ to the numbers $0,1,\dots, \lvert V(T) \rvert-1$ as in the order given by $\leq_{sc}$.
\end{definition}
\begin{example}\label{ex: sc order example}
We depict the sublevel-connected Morse labeling in the following example:
$$
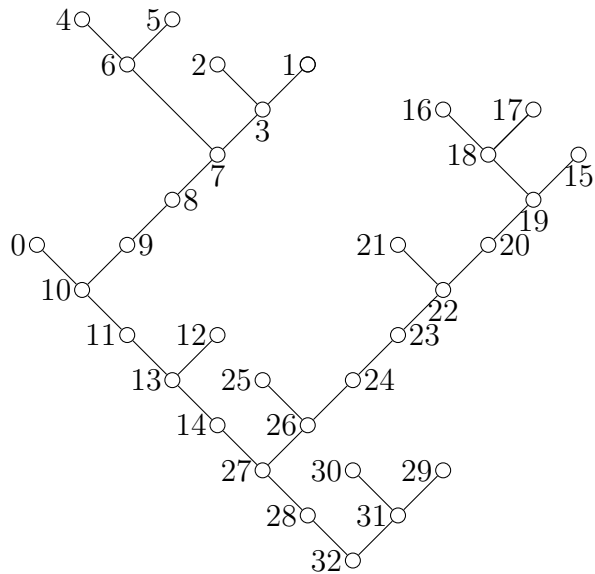
\begin{tikzpicture}[scale=0.60]

\node[inner sep=2pt, circle] (0) at (0,8) [draw] {};
\node[inner sep=2pt, circle] (8) at (1,7) [draw] {};
\node[inner sep=2pt, circle] (85) at (2,6) [draw] {};
\node[inner sep=2pt, circle] (10) at (3,5) [draw] {};
\node[inner sep=2pt, circle] (105) at (4,4) [draw] {};
\node[inner sep=2pt, circle] (20) at (5,3) [draw] {};
\node[inner sep=2pt, circle] (205) at (6,2) [draw] {};
\node[inner sep=2pt, circle] (24) at (7,1) [draw] {};
\node[inner sep=2pt, circle] (775) at (2,8) [draw] {};
\node[inner sep=2pt, circle] (75) at (3,9) [draw] {};
\node[inner sep=2pt, circle] (7) at (4,10) [draw] {};
\node[inner sep=2pt, circle] (6) at (5,11) [draw] {};
\node[inner sep=2pt, circle] (5) at (6,12) [draw] {};
\node[inner sep=2pt, circle] (3) at (2,12) [draw] {};
\node[inner sep=2pt, circle] (1) at (1,13) [draw] {};
\node[inner sep=2pt, circle] (2) at (3,13) [draw] {};
\node[inner sep=2pt, circle] (4) at (4,12) [draw] {};
\node[inner sep=2pt, circle] (5) at (6,12) [draw] {};
\node[inner sep=2pt, circle] (9) at (4,6) [draw] {};
\node[inner sep=2pt, circle] (19) at (6,4) [draw] {};
\node[inner sep=2pt, circle] (1775) at (7,5) [draw] {};
\node[inner sep=2pt, circle] (175) at (8,6) [draw] {};
\node[inner sep=2pt, circle] (17) at (9,7) [draw] {};
\node[inner sep=2pt, circle] (155) at (10,8) [draw] {};
\node[inner sep=2pt, circle] (15) at (11,9) [draw] {};
\node[inner sep=2pt, circle] (14) at (12,10) [draw] {};
\node[inner sep=2pt, circle] (18) at (5,5) [draw] {};
\node[inner sep=2pt, circle] (16) at (8,8) [draw] {};
\node[inner sep=2pt, circle] (13) at (10,10) [draw] {};
\node[inner sep=2pt, circle] (11) at (9,11) [draw] {};
\node[inner sep=2pt, circle] (12) at (11,11) [draw] {};
\node[inner sep=2pt, circle] (23) at (8,2) [draw] {};
\node[inner sep=2pt, circle] (22) at (9,3) [draw] {};
\node[inner sep=2pt, circle] (21) at (7,3) [draw] {};

\draw[-]  (0)--(8) node[midway, below] {};
\draw[-]  (775)--(8) node[midway, below] {};
\draw[-]  (85)--(8) node[midway, below] {};
\draw[-]  (775)--(75) node[midway, below] {};
\draw[-]  (75)--(7) node[midway, below] {};
\draw[-]  (7)--(6) node[midway, below] {};
\draw[-]  (7)--(3) node[midway, below] {};
\draw[-]  (2)--(3) node[midway, below] {};
\draw[-]  (6)--(4) node[midway, below] {};
\draw[-]  (6)--(5) node[midway, below] {};
\draw[-]  (1)--(3) node[midway, below] {};
\draw[-]  (85)--(10) node[midway, below] {};
\draw[-]  (10)--(9) node[midway, below] {};
\draw[-]  (10)--(105) node[midway, below] {};
\draw[-]  (105)--(20) node[midway, below] {};
\draw[-]  (20)--(19) node[midway, below] {};
\draw[-]  (19)--(18) node[midway, below] {};
\draw[-]  (19)--(1775) node[midway, below] {};
\draw[-]  (1775)--(175) node[midway, below] {};
\draw[-]  (175)--(17) node[midway, below] {};
\draw[-]  (17)--(16) node[midway, below] {};
\draw[-]  (17)--(155) node[midway, below] {};
\draw[-]  (155)--(15) node[midway, below] {};
\draw[-]  (15)--(14) node[midway, below] {};
\draw[-]  (15)--(13) node[midway, below] {};
\draw[-]  (13)--(11) node[midway, below] {};
\draw[-]  (13)--(12) node[midway, below] {};
\draw[-]  (20)--(205) node[midway, below] {};
\draw[-]  (205)--(24) node[midway, below] {};
\draw[-]  (24)--(23) node[midway, below] {};
\draw[-]  (13)--(12) node[midway, below] {};
\draw[-]  (23)--(21) node[midway, below] {};
\draw[-]  (23)--(22) node[midway, below] {};

\node[anchor = east]  at (0) {{$0L$}};
\node[anchor = east]  at (1) {{$4L$}};
\node[anchor = east]  at (2) {{$5R$}};
\node[anchor = east]  at (3) {{$6L$}};
\node[anchor = east]  at (4) {{$2L$}};
\node[anchor = east]  at (5) {{$1R$}};
\node[anchor = west]  at (6) {{$3R$}};
\node[anchor = west]  at (7) {{$7R$}};
\node[anchor = east]  at (8) {{$10L$}};
\node[anchor = east]  at (9) {{$12R$}};
\node[anchor = east]  at (10) {{$13L$}};
\node[anchor = east]  at (11) {{$16L$}};
\node[anchor = east]  at (12) {{$17R$}};
\node[anchor = east]  at (13) {{$18L$}};
\node[anchor = west]  at (14) {{$15R$}};
\node[anchor = west]  at (15) {{$19R$}};
\node[anchor = east]  at (16) {{$21L$}};
\node[anchor = west]  at (17) {{$22R$}};
\node[anchor = east]  at (18) {{$25L$}};
\node[anchor = east]  at (19) {{$26R$}};
\node[anchor = east]  at (20) {{$27L$}};
\node[anchor = east]  at (21) {{$30L$}};
\node[anchor = east]  at (22) {{$29R$}};
\node[anchor = east]  at (23) {{$31R$}};
\node[anchor = east]  at (24) {{$32L$}};

\node[anchor = west]  at (75) {{$8R$}};
\node[anchor = west]  at (775) {{$9R$}};
\node[anchor = east]  at (85) {{$11L$}};
\node[anchor = east]  at (105) {{$14L$}};
\node[anchor = west]  at (155) {{$20R$}};
\node[anchor = west]  at (175) {{$23R$}};
\node[anchor = west]  at (1775) {{$24R$}};
\node[anchor = east]  at (205) {{$28L$}};

\end{tikzpicture}
$$
We consider some nodes and their corresponding path words:
$$
\begin{tabular}{c|c}
   Label  & Path Word \\ \hline
    0 & LLLLLLLL \underbar{\mbox{ \ }} \underbar{\mbox{ \ }} \underbar{\mbox{ \ }} \underbar{\mbox{ \ }} \\
    1 & LLLLLLLRRRRR \\
    10 & LLLLLLL \underbar{\mbox{ \ }}  \underbar{\mbox{ \ }} \underbar{\mbox{ \ }}  \underbar{\mbox{ \ }}  \underbar{\mbox{ \ }} \\
    16 & LLLRRRRRRLL \underbar{\mbox{ \ }} \\
    32 & L \underbar{\mbox{ \ }} \underbar{\mbox{ \ }} \underbar{\mbox{ \ }} \underbar{\mbox{ \ }} \underbar{\mbox{ \ }} \underbar{\mbox{ \ }} \underbar{\mbox{ \ }} \underbar{\mbox{ \ }} \underbar{\mbox{ \ }} \underbar{\mbox{ \ }} \underbar{\mbox{ \ }} 
\end{tabular} 
$$
The node labeled 32 is maximal because it satisfies case b) with respect to all other nodes. By, \Cref{scmo}, the order relation between two arbitrary nodes is always decided at the level of their youngest common ancestor because the corresponding path words agree up to the position of the youngest common ancestor. 
\end{example}
\begin{proposition}\label{mergevsmo}
The construction of the sublevel-connected Morse labeling and forgetting the Morse labeling defines a pair of inverse bijections 
$$
\begin{tikzpicture}
\node (0) at (0,0) {$\lambda_{sc} \colon gMer/_{\cong}$}; 
\node (1) at (4,0) {$gMLT/_{\sim}\colon \text{forget}$};
\draw[<->] (0) -- (1);
\end{tikzpicture}
$$
where $\sim$ denotes merge equivalence and $gMer$ denotes isomorphism classes of generalized merge trees.
\end{proposition}
\begin{proof}
The statement follows directly by \Cref{corollary merge equivalence}.
\end{proof}

To summarize our results of this section, we take a look at how \Cref{DmfongraphvsMltree} and \Cref{mergevsmo} turn the different maps from \Cref{diagramminverseproblem} into bijections by dividing out the needed notion of equivalence. If we do not divide out any equivalence relation, the map $\Phi$ is not even well-defined. The maps $M(\_ \ ,\_)$ and forget are surjective, but not injective. The map $\lambda_{sc}$ is injective but not surjective. \par
Identifying CM-equivalent dMfs makes $\Phi$ a well-defined map and, moreover, a bijection which is inverse to $M(\_ \ ,\_)\colon DMF^{crit}_{graphs}\rightarrow gMLT$ by \Cref{DmfongraphvsMltree}. 
% Inverting order equivalences turns $iMO$ and $iML$ into inverse bijections. Finally,
Inverting merge equivalences makes $\lambda_{sc}$ and forget inverse to each other. As a consequence, we have a complete description of the inverse problem for critical dMfs on multigraphs and their induced merge trees. The characterization for arbitrary dMfs on 1-dim regular CW complexes follows by collapsing matched cells and then applying a version of \Cref{DmfongraphvsMltree} that incorporates \Cref{extend}. However, this procedure secretly makes use of two features which might become problematic if one tries to generalize the result to higher dimensions: on one hand, we use that irregularities of attaching maps can be easily characterized in the 1-dimensional case: here they always produce self-loops. Dealing with irregular faces in higher dimensions would be more difficult.\par
On the other hand, even if we start with regular CW complexes, the complex that arises by performing the simple collapses described by a Morse matching is not arbitrary but subject to being simple homotopy equivalent to a regular CW complex. It is a feature of dimension one that all 1-dimensional CW complexes are simple homotopy equivalent to a 1-dimensional regular CW complex. Hence, defining CM equivalences becomes more difficult in a higher-dimensional setting, in particular, if one wants to work with non-critical dMfs. This would lead to the need to analyze which CW complexes are simple homotopy equivalent to regular CW complexes in order to know for which generality a notion of CM equivalence is needed.

\section{Realization Problem with Simple Graphs}

Let $T$ be a generalized merge tree. Recall that $C(T)=C$ denotes the set of all cycle nodes of $T$.  For any $c\in C$, let $c_u$ denote the unique child of $c$.  For any $v\in T$, let $T(v)$ denote the generalized merge subtree\footnote{I.e. the subtree with root $v$ together with the chiralities inhertited from $T$.} of $T$ with root $v$  and let $\ell(v)$ denote the number of leafs of $T(v)$.

\begin{theorem}\label{thm: gmt chracterization}  Let $T$ be a generalized merge tree.  Then there exists a simple graph $X$ and  dMf $f\colon X\to \mathbb{R}$ such that $M(X,f)=T$ if and only if for every $c\in C(T)$, 

$$|C(T(c_u))|<\frac{(\ell(c_u)-2)(\ell(c_u)-1)}{2}.$$
Furthermore, $X$ can be made planar if and only if 
$$|C(T(c_u))|<2\cdot\ell(c_u)-5.$$
\end{theorem}

\begin{proof}
Suppose there exists a simple graph $X$ and dMf $f\colon X \to \mathbb{R}$ such that $M(X,f)=T$, and suppose by contradiction that there is a $c\in C(T)$ with the property that 
$$
|C(T(c_u))|\geq \frac{(\ell(c_u)-2)(\ell(c_u)-1)}{2}.
$$

By \Cref{lem: subtree component}, the rooted subtree $T(c_u)$ is isomorphic as rooted binary trees to the induced Morse labeled merge tree of $X_{f(s)}[s]$ where $s$ is the simplex of $X$ such that $M(s)=c_u$.  
Letting $v$ be the number of vertices in $X_{f(s)}[s]$, $e$ the number of edges in $X_{f(s)}[s]$, and $b_1$ the number of cycles in $X_{f(s)}[s]$, we see that \begin{eqnarray*}
e&=&v-1+b_1\\
&\geq&v-1+\frac{(v-1)(v-2)}{2}\\
&=&v-1+\frac{v(v-1)}{2}+1-v\\
&=& \frac{v(v-1)}{2}
\end{eqnarray*}
which is the maximum number of edges any connected component can have.  Hence it is impossible to add a cycle to this connected component so that $$
|C(T(c_u))|< \frac{(\ell(c_u)-2)(\ell(c_u)-1)}{2}.
$$
for all $c\in C$.  Now suppose further that $X$ is planar, and suppose by contradiction that $|C(T(c_u))|\geq 2\cdot\ell(c_u)-5.$ Using the same notation as above, we have 
\begin{eqnarray*}
e&=&v-1+b_1\\
&\geq&v-1+2v-5\\
&=& 3v-6.
\end{eqnarray*}
But it is well known that a simple planar graph satisfies $e\leq 3v-6$ \cite[Theorem 5.9]{Bickle2020}.  Hence either $X_{f(s)}[s]$ is not planar or maximal planar in the case of equality.  In either case, another edge cannot be added to $X_{f(s)}[s]$ without breaking planarity, and thus the result.
\bigskip

For the other direction, given the generalized Merge tree $T$, construct the sublevel-connected Morse labeling $\lambda_{sc}$ (\Cref{scmo}) on the nodes of $T$. 
% Associate to this Morse order a Morse labeling $\lambda \colon T\to \mathbb{R}$ such that $a\leq_{sc} b$ if and only if $\lambda(a)\leq \lambda(b)$. 
Consider the underlying Morse labeled merge tree $(\overline{T}, \overline{\lambda_{sc}})$ of $(T,\lambda_{sc})$, see \Cref{Def:CycleNodesandUnderlyingMLTree}.
% Apply the construction in \Cref{idmf} to $(T,\lambda)$ to obtain the underlying merge tree $(\overline{T}, \overline{\lambda})$.  
By \cite[Theorem 6.5]{brueggemann2021merge}, there is a path $P$ and dMf $\overline{f}\colon P\to \mathbb{R}$ such that $M(P,\overline{f})=(\overline{T}, \overline{\lambda_{sc}})$. We will inductively attach edges to $P$ in one-to-one correspondence with cycle nodes of $T$.  Each edge will be labeled with the same label as its corresponding cycle node. 

Induce on the cycle nodes of $T$ with respect to the sublevel-connected Morse labeling $c_1\leq_{sc}c_2\leq_{sc}\cdots$. For the base case $i=1$, write $P=X^1$.  We have by hypothesis that 
$$
|C(T({c_1}_u))|< \frac{(\ell({c_1}_u)-2)(\ell({c_1}_u-1)}{2}.
$$
In addition, $M(P,\overline{f})=(\overline{T}, \overline{\lambda})$ so ${c_1}_u=M(s_1)$ for some simplex $s_1\in P=X^1$.  Applying the correspondence noted in \Cref{rem graph and tree correspondence}, this inequality means that 
$$
b_1(X^1[s_1])< \frac{(v(X^1[s_1]-2)(v(X^1[s_1]-1)}{2}.
$$
By the computation in the forward direction, this implies that $e(X^1[s_1])<\frac{v(X^1[s_1])(v(X^1)-1)}{2}$.  Hence there are at least two vertices in $X^1[s_1]$ not connected by an edge.  A choice of vertex can be made by defining a lexicographic ordering on a subset of ordered pairs of the vertex set of $P$ where an ordered pair $(v,u)$ satisfies $\overline{f}(v)<\overline{f}(u)$ and $(v,u)<(v',u')$ if $\overline{f}(v)<\overline{f}(v')$ or $\overline{f}(u)<\overline{f}(u')$ when $\overline{f}(v)=\overline{f}(v')$.  Since all the vertices of $P$ are given distinct values, $<$ is a total order.  Add an edge $e_1$ incident with the vertices in the minimum pair over all available pairs to create $X^2=X^1\cup \{e_1\}$ and extend $\overline{f}$ to $f^1(e_1):=\lambda(c_1)$. Then $M(X^2,f^1)\simeq (T_{\leq \lambda(c_1)}, \lambda|_{T_{\leq \lambda(c_1)}})$. The inductive step is identical to the base case.

Now suppose that $|C(T(c_u))|<2\cdot\ell(c_u)-5$ for all cycle nodes $c\in T$.  By the forward direction, this is equivalent to $e<3v-6$ in the corresponding sublevel complex of $X$. The method of construction is analogous to the above construction and utilizes the fact that if a planar simple graph satisfies $e<3v-6$, then it is not maximal planar and hence an edge can be added while maintaining planarity \cite[Corollary 5.11]{Bickle2020}.
\end{proof}

\begin{remark} While the choices made in the construction of the simple graph $X$ in \Cref{thm: gmt chracterization} may be thought of as one canonical choice, the sublevel-connected Morse labeling is only one possible representative for the Morse labeling. Another just as natural (and shuffle equivalent\footnote{That is, Morse labelings that induce the same restricted order on leafs as well as the same restricted order on inner nodes. See \cite[Definition 2.24]{brueggemann2021merge} for details.}) labeling would be the index Morse labeling \cite[Definition 3.3]{brueggemann2021merge}. Furthermore, once a Morse labeling is picked, there are often several possible simple graphs with dMfs all related by CM equivalence that represent the given generalized merge tree. 
    
\end{remark}

\begin{example} To illustrate the construction in the planar case, consider the generalized merge tree $T$ pictured below:

$$

$$

We constructed the sublevel-connected Morse labeling $\lambda_{sc}$ in \Cref{ex: sc order example}.

$$
%
$$

We then pass to the underlying merge tree $\overline{T}$ and restrict $\lambda$ to $\overline{T}$ in order to apply \cite[Definition 4.5]{brueggemann2021merge} to obtain the sublevel-connected dMf on the graph below with induced merge tree $\overline{T}$.

$$
%
$$

We induce on the cycle nodes ordered by their generalized Morse label. The first cycle to be introduced is cycle node with label 8. This will be a cycle added to the graph

$$
%
$$
to the component with the edge labeled $7$. 

$$
%
$$
We then add the cycle corresponding to the node labeled 9  to this same graph.
$$
%
$$

Skipping to the cycle node labeled 23, we see that we need to add a cycle to the component with edge labeled 22:

$$
%
$$
We add this edge

$$
%
$$
and must add another cycle corresponding to cycle node labeled 24 to this same connected component.

$$
%
$$

Notice that this component is now a complete graph and that no more cycles can be added.  The final graph with dMf that induces the given generalized merge tree is
$$
%
$$

\end{example}

\section{How to Find Cancellations with Merge Trees}\label{section cancel critical}
One of the desirable features of discrete Morse functions is that they simplify the computation of cellular homology: there is a chain complex, called the Morse complex, that has only the critical cells of a discrete Morse function as generators and is chain-equivalent to the cellular chain complex \cite[Section 7]{Forman98}. Thus discrete Morse functions with as few critical cells as possible are useful for efficient computation of ordinary homology of cell complexes. Finding optimal discrete Morse functions, i.e., ones which minimize the number of critical cells, is NP-complete \cite{HaJa23}. A common approach to find ``good" discrete Morse functions, i.e., ones which have relatively few critical cells but are not necessarily optimal, is to start with any discrete Morse function and to modify it by canceling critical cells. Canceling critical cells refers to the process of inverting unique gradient paths between pairs of critical cells, see \cite[Theorem 11.1]{Forman98}, which leads to a new combinatorial gradient field for which said pairs of cells are no longer critical. Each combinatorial gradient field is represented by discrete Morse functions, so we can also find a new discrete Morse function that realizes the matching. \par
Before we continue with the algorithm, we recall the concepts of combinatorial gradient fields and gradient paths.
\begin{definition}[{\cite[Section 3, Definition 3.3]{Forman2002}}]
    Let $X$ be a graph. A \emph{discrete vector field} $V$ on $X$ is a collection of pairs $V\subset X \times X$ such that if $(\sigma,\tau) \in V$, then $\sigma \subsetneq \tau$ and each cell is contained in at most one pair of $V$. \par
    Let $f\colon X \rightarrow \mathbb{R}$ be a discrete Morse function on $X$. The \emph{combinatorial gradient field} $\nabla f$ of $f$ is the discrete vector field that consists of all the pairs $(\sigma,\tau)$ such that $f(\sigma)=f(\tau)$.\par
    For any discrete vector field $V$ on $X$, a \emph{$V$-path} is a sequence of cells 
    $$
    \alpha_0, \beta_0, \alpha_1, \beta_1, \dots \beta _r,\alpha_{r+1}
    $$
    such that $(\alpha_i,\beta_i) \in V$ and $\alpha_i\neq \alpha_{i+1} \subset \beta_{i}$ for all $i$. If $V$ is the combinatorial gradient field of a discrete Morse function $f$, then we call the $V$-paths the \emph{gradient paths} of $f$.
\end{definition}
In this section, we present a way to find cancellations of critical cells of dMfs with the help of the induced ML tree while preserving certain potentially relevant features. This presentation is meant to be a purely theoretical contribution for this work because in the case of graphs, similar results can be obtained by usage of efficient algorithms that find spanning trees of the given graph together with techniques to find optimal dMfs on trees, see e.g.~\cite{RandScoville2020}. We conjecture that the following ideas might lead to useful simplifications of dMfs on higher dimensional complexes, once a generalization of \Cref{iMLT} to higher dimensional complexes will be found. \par
 
% In order to find an optimal dMf on any given graph, one may start with an arbitrary dMf that only has critical cells and continue with performing cancellations along the merge tree. 
The main idea of the algorithm is to perform cancellations of critical cells along the induced ML tree of a given dMf. In case no dMf is previously given, one may choose any arbitrary dMf for which every cell is critical. 
\begin{remark}
In order to obtain an arbitrary dMf on a graph $X$ that has only critical cells, one can simply choose any total order on the vertices and any total order on the edges. Then assign the values $0,\dots, \lvert V(X) \rvert -1$ to the vertices according to the chosen order and the numbers $\lvert V(X) \rvert, \dots, \lvert V(X) \rvert+ \lvert E(X) \rvert $ to the edges. This always produces an index-ordered dMf, which is not necessary for the following algorithm. Perhaps more sophisticated approaches to finding a critical dMf might be useful, but for now we are satisfied with this simple one.
\end{remark}
The input of the algorithm is a dMf $f\colon X \rightarrow \mathbb{R}$ on a graph $X$. The output of the algorithm depends on user choices, but will always be a dMf $f'\colon X' \rightarrow \mathbb{R}$ on a graph $X'$ such that $X'$ is homotopy-equivalent to $X$, has the same number of vertices and edges as $X$, and $f'$ has a lower or equal number of critical cells compared to $f$.
Given a critical dMf $f\colon X \rightarrow \mathbb{R}$, the algorithm proceeds as follows:
\begin{enumerate}
    \item Calculate the induced generalized Morse labeled merge tree $M(X,f)$, and let $C$ be the set of leaves of $M(X,f)$
    \item If $C=\emptyset$, go to (3). Otherwise, let $c\in C$ be the leaf with maximal label. If $c$ has no ancestor that is neither a cycle node nor matched, change $C=C-\{c\}$ and return to (2).\footnote{This step is needed for the case that $c$ is the last remaining leaf.} Otherwise, let $p$ be the youngest ancestor of $c$ such that $p$ is neither a cycle node nor matched. Then either:
    \begin{enumerate}[a)]
        \item The vertex $M^{-1}(c)$ is adjacent to the edge\footnote{Recall that ancestors need to be inner nodes and, hence, correspond to edges of $X$. See \Cref{Rem:CorrespondenceBetweenNodesAndVerticesAndEdges}.} $M^{-1}(p)$. 
        \item The vertex $M^{-1}(c)$ is not adjacent to the edge $M^{-1}(p)$.
    \end{enumerate}
    If case a), match $M^{-1}(c)$ and $M^{-1}(p)$.
    This does not produce cycles because we explicitly exclude cycle nodes from the matching. Let $C=C-\{c\}$ and return to (2) \\
   
    If case b), either:
    \begin{enumerate}[i)]
        \item leave $M^{-1}(c)$ critical, let $C=C-\{c\}$ and return to (2)
        \item check for a symmetry equivalence $a$ of $(X,f)$ such that $a(M^{-1}(c))$ is adjacent to $a(M^{-1}(p))$, apply it, and then proceed as in case a).  If there is no symmetry equivalence, proceed to i), iii), or iv).
        \item apply a CM equivalence in order to make $M^{-1}(c)$ and $M^{-1}(p)$ adjacent, then proceed as in case a), or
        \item observe that there is a unique gradient path from $M^{-1}(c)$ to $M^{-1}(p)$ and cancel the two cells along this gradient path. Let $C=C-\{c\}$ and return to (2).
    \end{enumerate}
    \item At this step we have a combinatorial gradient field $V$ on a graph $X'$. In order to define $f'$, we assign each critical cell $c$ of $V$ the value $f(c)$. We further extend $f'$ according to $V$ as follows: let $c$ be the smallest critical cell such that gradient paths $\gamma$ end at $c$ for which not all matched cells have been assigned values. Let $c'$ be the smallest critical cell such that one of these gradient paths to $c$ begins at $c'$. Then we assign to all pairs of matched cells on the path from $c'$ to $c$ values between $f'(c)$ and $f'(c')$ in descending order. For gradient paths that only end in at $c$ but do not start at a critical cell, we assign values in increasing order away from $c$.\footnote{The matched cells are part of gradient paths between the critical cells. Therefore, their values under $f'$ need to lie between the values of the critical cells at the beginning and end of the corresponding gradient path.} We perform this inductively over all critical cells. Then $f'$ is delivered as output.
\end{enumerate}
The precise nature of the output depends on the choices the user makes in case b).  If case b) never applies, the output will be an optimal discrete Morse function on the exact same graph $X$. In the case that b) is applied and the user chooses option i), an optimal matching is not guaranteed but we still preserve $X$. Case ii) is only available if a suitable symmetry equivalence actually exists, which does not necessarily need to be the case. Moreover, the existence of such symmetry equivalences needs to be checked for each level individually because automorphism groups of graphs can change arbitrarily along filtrations. Nonetheless, if it is possible to apply case ii) consistently, then the output discrete Morse function will be optimal and the graph will be unchanged. If case iii) is consistently chosen, we produce an optimal matching but may change the homeomorphism type of $X$. If case iv) is consistently chosen, we preserve $X$ and obtain an optimal matching but we change the order of the vertices induced by $f$ on a larger scale. While one could in principle choose different options of case b) at different stages in a single run of the algorithm, this would produce a seemingly undesirable output, as it would suffer all the drawbacks mentioned in each case.  

Most of the claims made in the above algorithm are straightforward to prove. For example, the fact that the cases 2a), 2b)i), and 2b)iii) work as described follows immediately from the definition of the used equivalences. However in general it does not appear easy to decide whether case 2b)ii) is applicable. Nonetheless, case 2b)iv) is not so obvious, so we consider it in the following lemma:
\begin{lemma}
Let $X$ be a graph, $f\colon X\rightarrow \mathbb{R}$ a critical dMf, and $M(X,f)$ the induced gML tree. At any point of the cancellation algorithm, there is always a unique gradient path from the vertex $M^{-1}(c)$ corresponding to the maximally labeled unmatched leaf $c$ to the edge $M^{-1}(c)$ corresponding to its youngest unmatched ancestor $p$.  
\end{lemma}
\begin{proof}
If $M^{-1}(c)$ and $M^{-1}(p)$ are adjacent, there is nothing to prove. If $M^{-1}(c)$ and $M^{-1}(p)$ are not adjacent then there is no non-closing critical edge in $X_{f(M^{-1}(p)-\varepsilon)}[M^{-1}(c)]$ because otherwise said other younger critical edge would induce a younger unmatched ancestor of $c$. \par
Since $M^{-1}(c)$ is a critical vertex with no adjacent critical edge, all adjacent edges of $M^{-1}(c)$ are matched with their respective other vertex. This means that on all adjacent edges, there is a gradient path pointing towards $M^{-1}(c)$. Following these gradient paths backwards either leads to matched vertices that are adjacent only to the edge they are matched with, or to the unique non-closing critical edge of $X_{f(M^{-1}(p))}[M^{-1}(c)]$. One of the gradient paths eventually leads to $M^{-1}(p)$ because $X_{f(M^{-1}(p))}[M^{-1}(c)]$ is connected.

The gradient path is unique because closing edges remain critical, that is, because we only match cells along a subtree of $X$.
\end{proof}

We apply the cancellation algorithm to the following example:
\begin{example}
We consider the graph:\\
$$

$$
We put some critical dMf on it and calculate the induced generalized merge tree:
$$
%
$$
We perform the algorithm as long as only step 2a) occurs:

$$
%
$$
Now is the first time we run into case 2b). We can actually apply case 2b)ii) here:

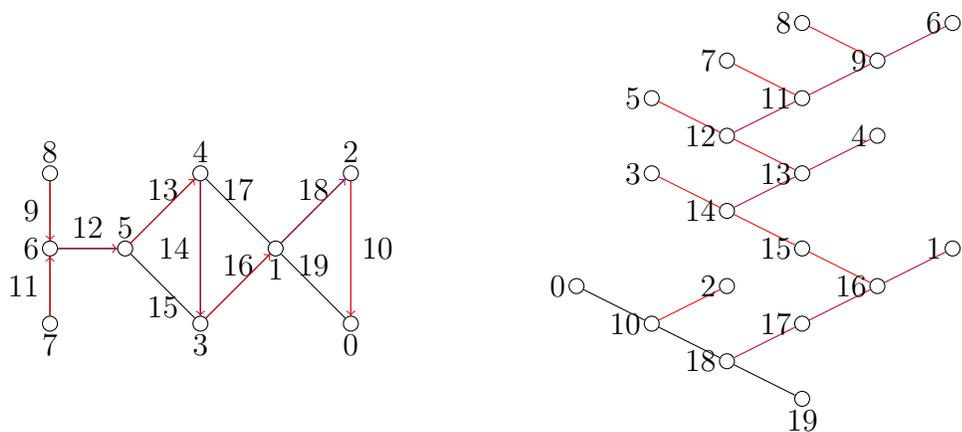
$$
%
$$
At this point we have finished the construction of the new combinatorial gradient field $V$. A possible discrete Morse function $f'$ for $V$ is:
$$
%
$$
In this example, the cases 2a) and 2b)ii) sufficed.

\end{example}
We consider the following example in order to see how quickly things can fail:
\begin{example}
We consider the following dMf and its induced merge tree:
$$
%
$$
The algorithm runs into case 2a) four times, which results in the following:
$$
%
$$
Now we have reached case 2b) and case 2b)ii) is not applicable. We would need to have the vertex labeled 1 adjacent to the edge labeled 10. But this is not possible because all symmetry equivalences leave the vertex labeled 1 adjacent to the edge labeled 9 and no other edge. The three different options lead to the following:
$$
%
$$
Possible choices for representing discrete Morse functions $f'$ could be:
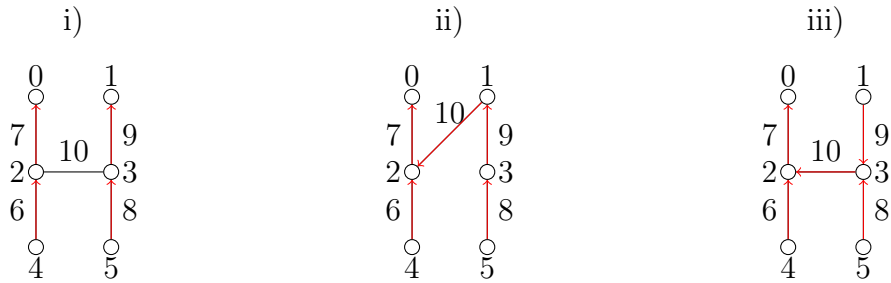
$$
%
$$
\end{example}
\begin{example}
A sublevel symmetry of the last sublevel complex before the ``merge tree algorithm" fails may not always be sufficient. Consider the graph with dMf given below. 

$$
%
$$

Proceeding as before, we obtain a matching on the graph until the algorithm specifies to match the vertex labeled $1$ with the edge labeled $13$.  Since these cells are not incident, we need to find a sublevel-symmetry of sublevel $12$.  However, the sublevel subcomplex $X_{12}$ is given by
$$
%
$$
which is well-known to have no non-trivial automorphisms. There is also no symmetry equivalence of a lower level than 12 that makes the vertex labeled 1 and the edge labeled 13 adjacent. However, the three different workarounds mentioned earlier result in the following:

$$
%
$$

\end{example}

At the end of this section, we compare our algorithm for finding cancellations of critical cells to similar algorithms from the literature. In \cite{LewinerLopesTavares2003} the authors introduce an algorithm to find optimal dMfs on 2-dimensional manifolds which they generalize to higher dimensions and more general complexes in \cite{LEWINER2003221}, even though losing the guarantee for optimality in the process. The main similarity to our approach is the use of an auxiliary tree structure, in our case the generalized merge tree, in the case of \cite{LEWINER2003221} a spanning hyperforest of a hypergraph associated to the Hasse diagram of a dMf.  \par
In \cite{RandScoville2020}, the authors provide an algorithm to find optimal dMfs on trees. Said algorithm, combined with any standard algorithm to find spanning trees, can easily be generalized to provide optimal dMfs on graphs with a prescribed critical vertex.\par
The main feature of our new approach, compared to the pre-existing ones, seems to be that our algorithm allows to preserve certain properties of a given dMf. In certain cases, such a dMf might be given by an application and, therefore, might be worth preserving. We conjecture that, given a suitable version of higher merge trees, our algorithm can be generalized to higher dimensions. Since finding optimal Morse matchings is MAX–SNP hard, such a generalization might either fail to be optimal or be inconvenient to work with in practice. Nonetheless, we hope to find interesting classes of examples in which such a generalized algorithm happens to be performative and informative. 

\section{Future Directions}
In this section, we want to take a look at possible applications and further directions this work might lead to.\par
Our main results, \Cref{DmfongraphvsMltree} and \Cref{thm: gmt chracterization}, give a detailed description of the fiber of the persistence map that takes dMfs on graphs to their persistent connectivity. This approach may be used in applications in which the persistent connectivity turns out to be the most relevant feature, allowing to replace a maybe inconvenient graph with a more convenient one that describes that same persistent connectivity. At the end of \Cref{section preliminaries} we sketch how to search through all possible representatives in a structured way. \Cref{thm: gmt chracterization} provides an easy-to-check condition for when this replacement can be chosen to be a simple graph. \par
A similar approach is given by applying the cancellation algorithm from \Cref{section cancel critical}. The algorithm helps to simplify dMfs on graphs while allowing to preserve either the homeomorphism type or the dynamics induced by the Morse function. 
One immediate question would be, how much approaches such as these change the original Morse function. Thus, it seems interesting to investigate the diameter of the set of representatives for a given merge class of dMfs with respect to some suitable metric for dMfs. Moreover, we would be interested in finding out how distant the function coming from the cancellation algorithm is from its input function.\par
In a more pure direction, one could try to set up a version of persistent geometric group theory using the groups of sublevel automorphisms as in \Cref{examplesofslsymm}. On one hand, it seems interesting in itself to consider actions of sequences of groups on sequences of spaces and which ones can be realized as sublevel automorphisms of a filtered space. On the other hand, the results from such approaches would be useful for applications of the cancellation algorithm mentioned above. It also seems interesting to analyze how the application of CM equivalences affects the sublevel automorphism group.\par
Furthermore, the notions of symmetry equivalences and CM equivalences might be helpful for the investigation of the space of dMfs on some given graph. \par
The most straightforward direction would be a generalization of \Cref{iMLT} to higher dimensions in order to enable the pursuit of all the above mentioned possible future directions in higher dimensions.

\vspace{5mm}

\section{Conflict of Interest Statement}
	The authors state that there is no conflict of interest.
\section{Data Availability Statement}
Data sharing is not applicable to this article as no datasets were generated or analyzed during the current study.
\bibliographystyle{alpha}
\bibliography{References}
\end{document}